\documentclass[11pt]{article}
\usepackage[a4paper,margin=1in]{geometry}
\usepackage{amsmath, amssymb, amsthm}
\usepackage{algorithm}
\usepackage{algpseudocode}
\usepackage{setspace}
\usepackage{hyperref}
\usepackage{graphicx}
\usepackage{subcaption}
\usepackage{enumitem}
\usepackage[svgnames]{xcolor}
\usepackage{booktabs}
\usepackage{caption}
\usepackage{placeins}

\newtheorem{theorem}{Theorem}[section]
\newtheorem{definition}[theorem]{Definition}

\newtheorem{proposition}[theorem]{Proposition}
\newtheorem{corollary}[theorem]{Corollary}

\newtheorem{assumption}[theorem]{Assumption}

\newcommand{\grad}{\nabla}

\newcommand{\R}{\mathbb{R}}

\DeclareMathOperator*{\argmin}{arg\,min} 

\newcommand{\revision}[1]{{\color{black}#1}}
\newcommand{\liurevision}[1]{{\color{black}#1}}
\title{Robust Estimators for Stochastic Model-Based Derivative-Free Optimization}
\author{
Fangyu Liu\thanks{
Corresponding author.
School of Mathematics and Statistics,
University of Melbourne,
Parkville VIC 3010, Australia.
ORCID: 0009-0004-1216-7307
(\texttt{fangyu.liux@gmail.com}).
This author was supported by the Australian Research Council
Discovery Early Career Award DE240100006.
}
\and
Lindon Roberts\thanks{
School of Mathematics and Statistics
\& ARC Training Centre in Optimisation Technologies,
Integrated Methodologies, and Applications (OPTIMA),
University of Melbourne,
Parkville VIC 3010, Australia.
ORCID: 0000-0001-6438-9703
(\texttt{lindon.roberts@unimelb.edu.au}).
This author was supported by the Australian Research Council
Discovery Early Career Award DE240100006.
}
}

\date{}
\begin{document}
\maketitle

\begin{abstract}
Derivative-free optimization (DFO) methods are popular for solving stochastic optimization problems, with model-based methods being a practically successful class of DFO algorithm. 
In most existing works, sufficiently accurate objective estimates (used for model construction and to test for sufficient decrease) are constructed by taking a sample mean of several noisy evaluations. 
In this work, we demonstrate that using robust estimators (specifically the median, trimmed mean, and median-of-means) can improve the sample complexity of model-based DFO methods. Our numerical results demonstrate that robust estimators are also a practical choice for algorithm implementations.

\end{abstract}

\medskip



\section{Introduction}
\label{sec:introduction}

Derivative-free optimization (DFO) refers to optimization methods that
use only function values
\cite{conn2009introduction,larson2019derivative}.
These methods are useful when derivative information is unavailable,
unreliable, or too expensive to compute.
This situation arises in many black-box, simulation-based, and noisy
optimization problems.
Model-based DFO \cite{roberts2025interpolation} is a widely used class
of DFO methods.
At each iteration, a local model of the objective is constructed from
sampled function values.
The model is then minimized inside a trust region.
The trust region controls the region in which the local model is
expected to be reliable and provides a framework for global convergence
\cite{conn2009introduction}.

Many applications involve stochastic objective evaluations.
In simulation optimization, for example, the problem often has the form
\begin{align}
    \min_{x\in\R^d} f(x)
    =
    \mathbb{E}_{\xi}[F(x,\xi)].
    \label{eq_simopt}
\end{align}
The algorithm cannot evaluate \(f(x)\) exactly.
Instead, it observes noisy realizations \(F(x,\xi)\).
Examples of applications involving stochastic optimization or noisy
evaluations include simulation-based inventory management
\cite{ye2016inventory}, climate model calibration
\cite{howland2022parameter}, and variational quantum optimization
with measurement noise \cite{scriva2024challenges}.

Noisy or unreliable objective values can also arise from numerical and
hardware failures.
Such failures can be important when each simulation is expensive.
For example, the classical 31-variable helicopter rotor design problem
studied in \cite{booker1998helicopter} involved frequent failures of
objective evaluations during the simulation process.
Large computing systems can also experience hardware failures.
A 128-CPU Beowulf cluster at Los Alamos National Laboratory, for
example, was reported to fail approximately once per week under normal
winter operating temperatures \cite{feng2008green}.
In some applications, a failed computation returns an explicit error
flag.
The failure can then be treated as a hidden constraint
\cite{chen2016hidden,ledigabel2024taxonomy}.
In other cases, the computation may return an inaccurate value without
warning.
A code may also construct an approximate output from the information
that remains available.
Fault-tolerant sparse-grid combination methods provide one example of
this behavior \cite{harding2015fault}.
These types of failures are more difficult to handle because the
returned value may be far from the value that would have been obtained
from a successful evaluation, \revision{and the returned value may not follow an easy-to-describe distribution}.

Accurate objective estimates are important in model-based DFO for two
reasons.
First, sampled values are used to construct a local interpolation or
regression model.
Second, function estimates are used to decide whether a trial point
gives sufficient improvement over the current iterate.
In stochastic model-based DFO, repeated noisy evaluations at the same
point are often combined using sample averaging to produce a more
accurate objective estimate.
The resulting estimate is then used for model construction or step
acceptance.

In this work, we study robust alternatives to sample averaging.
We consider the sample median, trimmed mean, and median-of-means.
Under finite-variance noise, robust mean estimators such as the trimmed
mean and median-of-means can achieve high-confidence accuracy with
better dependence on the required confidence level than the Chebyshev
guarantee for the sample mean.
The sample median also has exponential concentration when the
distribution has sufficient probability mass around its population
median.
The difference is larger under some failure models.
In particular, the sample median can recover the noise-free function
value with high probability when successful evaluations form a
majority.
Our numerical experiments show that robust estimators perform
similarly to the sample mean under relatively benign noise, but can
perform better when extreme observations or evaluation failures occur.

\paragraph{Existing Works}

Several lines of work have studied model-based DFO with stochastic
information.
Bandeira, Scheinberg, and Vicente
\cite{bandeira2014convergence} established convergence when local
models are probabilistically accurate but the function values used to
assess sufficient decrease are deterministic.
Gratton, Royer, Vicente, and Zhang \cite{gratton2018complexity} later
derived worst-case complexity bounds for this setting.
Chen, Menickelly, and Scheinberg \cite{chen2018stochastic} extended the
framework to allow both random models and random function estimates.
Their STORM framework proves global convergence when the required
accuracy conditions hold with sufficiently high probability.
Blanchet, Cartis, Menickelly, and Scheinberg
\cite{blanchet2019stochastic} subsequently established expected
iteration-complexity bounds.
More recent work has considered broader probabilistic oracle
conditions \cite{cao2024first}.
For example, Scheinberg and Xie \cite{scheinberg2025unreliable} allow
less reliable stochastic information, including heavy-tailed function
estimates, and derive high-probability complexity guarantees.
Related extensions include composite nonsmooth optimization
\cite{baraldi2025proxstorm} and constrained optimization using SQP
methods \cite{fang2026heavy}.

The STORM framework itself does not prescribe how accurate function
estimates must be constructed.
Its standard stochastic DFO construction, however, uses sample
averages.
Under unbiased noise with uniformly bounded finite variance,
sufficiently accurate function estimates and interpolation models are
obtained using sample averaging together with Chebyshev's inequality
\cite{chen2018stochastic}.
The cost of constructing such stochastic information has also been
studied directly.
Jin, Scheinberg, and Xie \cite{jin2024samplecomplexity} analyze the
total oracle complexity of adaptive stochastic optimization methods,
including STORM.
Rinaldi, Vicente, and Zeffiro \cite{rinaldi2024stochastic} introduce
a weaker tail-bound condition on the stochastic reduction estimates
used in the sufficient-decrease test.
Their condition allows the number of samples required by stochastic
trust-region and direct-search methods to be reduced.
These works show that the cost of obtaining sufficiently reliable
stochastic information can have an important effect on the overall
algorithm.

Other model-based DFO methods handle noisy data through the
construction of the local model itself.
Billups, Larson, and Graf \cite{billups2013derivative} use weighted
regression when function evaluations have different levels of
computational error.
Larson and Billups \cite{larson2016stochastic} construct stochastic
trust-region models from noisy function values using regression-type
approximations.
Their method does not require repeated sampling at the same point.
Instead, local regression models are used both to construct the
trust-region model and to estimate the function values needed for the
step-acceptance test.
Under independent, mean-zero, finite-variance noise, they prove that
the gradients at the iterates converge in probability to zero.
The focus of these methods is therefore the construction of useful
local models from noisy observations.

A different approach is to adapt the number of samples during the
optimization.
Rather than fixing the number of replications in advance, adaptive
sampling methods increase the sampling effort when a more accurate
estimate is needed.
The variable-number sample-path method of Deng and Ferris
\cite{deng2009variable} is an early example in a trust-region setting.
The STRONG framework of Chang, Hong, and Wan \cite{chang2013strong}
combines trust-region ideas with response-surface methodology.
It uses simulation replications together with experimental design and
regression to construct local response-surface models.
The ASTRO-DF framework of Shashaani, Hashemi, and Pasupathy
\cite{shashaani2018astro} also chooses the Monte Carlo sample size
adaptively.
It balances sampling error against local model error and establishes
almost sure convergence to a first-order critical point.
Its original analysis assumes finite-variance noise together with an
additional bounded higher-moment condition.
More recent work has studied iteration and sample complexity
\cite{ha2025adaptive,ha2025complexity}, stratified adaptive sampling
\cite{amici2026stratified}, and adaptive regularization
\cite{ha2026adaptive}.
The complexity analyses in
\cite{ha2025adaptive,ha2025complexity} use sub-exponential-type
conditions on stochastic estimation errors.
By contrast, the stratified sampling analysis in
\cite{amici2026stratified} relaxes the need for sub-exponential
function-output noise under additional assumptions.
The adaptive regularization method in \cite{ha2026adaptive} uses both
stochastic function and gradient information and also imposes
sub-exponential-type noise conditions.
These methods mainly study how sampling effort and stochastic
information are managed as optimization progresses, \revision{and again are based on using the sample mean to approximate the mean of a (true) underlying distribution for a stochastic objective}.

Recent work has also studied the effect of persistent errors in
interpolation-based model-based DFO.
Chaudhry and Scheinberg \cite{chaudhry2025complexity} derive
function-evaluation complexity bounds for classical model-based
trust-region methods and several modern variants.
They also consider an inexact zeroth-order oracle with uniformly
bounded function-value error.
A fixed error level creates a corresponding lower limit on the
optimization accuracy that can be guaranteed.
In this setting, the method can only be expected to approach a
neighborhood of stationarity once the oracle error becomes dominant.
This differs from the STORM setting considered here, where the
required function-estimate error decreases with the trust-region
radius and can therefore support convergence to a stationary point.

Noisy DFO also appears in robust optimization and optimization under
uncertainty.
Menhorn et al.~\cite{menhorn2022trust} study constrained stochastic
DFO problems involving expectations, standard deviations, and
risk-related objective and constraint quantities.
Their method uses sample-based estimates, local surrogate models,
Gaussian-process smoothing, and noise-aware trust-region management.
Here, robustness mainly concerns the quantities being optimized and
the treatment of uncertainty, rather than replacing the sample mean
with a robust mean estimator.
The authors develop and analyze the estimation and smoothing
procedures, but leave a convergence analysis of the proposed
stochastic optimization framework for future work.

Robust aggregation has also appeared recently in zeroth-order
optimization.
Ma, Li, and Zhang \cite{ma2026projection} propose PRISM-ZO for
zeroth-order optimization on Riemannian manifolds under heavy-tailed
oracle noise.
Their method samples random low-dimensional subspaces of the tangent
space and uses directions from these subspaces to construct
zeroth-order gradient estimates.
It then applies median-of-means or Huber aggregation to these
estimates.
This provides a recent example in which robust aggregation is built
directly into a zeroth-order optimization method.
The setting is different from the one studied here.
PRISM-ZO aggregates gradient estimates in the tangent space, whereas
our focus is on aggregating repeated scalar function evaluations
before they are supplied to a model-based DFO algorithm.

Robust aggregation of scalar function estimates has also been studied
in stochastic constrained optimization.
Fang, Lavaei, and Na \cite{fang2026heavy} develop a stochastic
trust-region SQP method for problems with a stochastic objective and
deterministic equality constraints.
Their analysis uses STORM-type probabilistic accuracy conditions and
allows heavy-tailed noise.
In particular, they consider noise with a bounded
\((1+\epsilon)\)-th moment for \(\epsilon\in(0,1]\), which permits
infinite variance.
They analyze both sample averaging and median-of-means for
constructing accurate stochastic function estimates.
Their median-of-means bounds improve the dependence on the required
failure probability from polynomial to logarithmic.
The paper also establishes high-probability iteration-complexity
bounds.
Its numerical experiments, however, use sample averaging rather
than comparing different aggregation rules.

These studies use different approaches to handle noisy optimization
information.
Larson and Billups \cite{larson2016stochastic} use local regression,
while Ma, Li, and Zhang \cite{ma2026projection} apply robust aggregation
to zeroth-order gradient estimates.
Fang, Lavaei, and Na \cite{fang2026heavy} show how median-of-means can
improve the accuracy guarantees of stochastic function estimates in a
trust-region SQP framework.

We build on this perspective by studying the aggregation of repeated
scalar function evaluations in model-based DFO.
Our focus is on comparing several aggregation rules under different
noise assumptions, including failure noise.
We examine their sample-size requirements in relation to the
probabilistic accuracy conditions used in STORM-type methods.
We also investigate how the choice of aggregation rule affects the
performance of a model-based DFO solver.
This complements \revision{and significantly extends} the theoretical analysis \revision{from}
\cite{fang2026heavy}, whose numerical experiments do not directly
compare different aggregation rules.

\paragraph{Contributions}

We study how the choice of aggregation rule affects the cost of
producing accurate function estimates for stochastic model-based DFO.
We compare the sample mean with the sample median, trimmed mean, and
median-of-means.
Our aim is to examine these estimators within a common STORM-type
accuracy framework and to understand how their theoretical properties
affect the performance of a model-based DFO solver.

We first compare the sample-size requirements of the four estimators
under finite-variance noise.
The Chebyshev guarantee for the sample mean requires \revision{$N=O(\delta^{-2} (1-\beta)^{-1})$} samples to obtain a
\((\delta,\beta)\)-probabilistically accurate estimate in the sense of
Definition~\ref{def:prob-accurate}.
By contrast, the trimmed mean and median-of-means require $N=O(\delta^{-2} \log((1-\beta)^{-1}))$ samples, up to problem-dependent constants.
The sample median also has logarithmic dependence on
\revision{$(1-\beta)^{-1}$} when the distribution has sufficient probability mass
around its median.
Our comparison makes the implications of these established
concentration results for the accuracy requirements of stochastic
model-based DFO explicit.
Under the corresponding assumptions, the resulting guarantees for the
robust estimators have exponential dependence on sample size, whereas
the Chebyshev guarantee for the sample mean gives polynomial decay.

This difference can become particularly important when the local
model uses many interpolation points.
Suppose that \(p\) independent function estimates are used to
construct a model that must be accurate with probability at least
\(\alpha_m\).
It is sufficient in our setting for each estimate to have accuracy
probability \revision{$\beta=\alpha_m^{1/p}$.} As \(p\) increases, the required value of \(\beta\) becomes
increasingly close to one.
The dependence on \(1-\beta\) can therefore become important for
models with larger interpolation sets, including those used in higher
dimensions.

Our second contribution concerns failure noise.
Suppose a function evaluation returns the noise-free value with
probability \(q>1/2\).
The sample median can then recover the noise-free function value
exactly with exponentially high probability.
Its guarantee does not depend on the magnitude of the failed
evaluations.
A suitably chosen trimmed mean also gives exact recovery under the
fixed-outlier failure model.
The sample mean and median-of-means behave differently.
Their accuracy bounds depend on the magnitude of the failed value,
which may be arbitrarily large.

We then compare the four estimators numerically on the low-dimensional
Mor\'e--Wild problems and a selected collection of medium-dimensional
CUTEst problems.
Under Gaussian noise, the four estimators perform similarly.
Under finite-variance Student-\(t\) noise with three degrees of
freedom, the differences are generally small, although robust
aggregation gives some improvement on the selected CUTEst problems.
Under Student-\(t\) noise with one degree of freedom, the robust
estimators perform better than the sample mean, particularly on
the CUTEst problems.
Under failure noise, the sample median performs substantially better
when a fixed number of repeated evaluations is used at each point.
However, its advantage can be lost when the sample size is determined
by the trust-region radius.
These results show that the performance of an aggregation rule also
depends on how samples are allocated during optimization.

Based on our analysis and numerical experiments, we make several
recommendations on the choice of aggregation rule.
The appropriate estimator depends on the population quantity that the
optimization problem is intended to minimize.
When the population median is consistent with the desired objective
in the sense of Assumption~\ref{ass_consistency}, the sample median is
a strong robust choice, especially when failures or extreme
observations are possible.
If the problem is intrinsically defined by a mean, as in
\eqref{eq_simopt}, changing from the sample mean to the sample median
may change the optimization objective.
In this case, median-of-means provides a natural robust alternative
under finite-variance noise.
A suitably constructed trimmed mean provides another option when its
bias can be controlled.
Under sub-Gaussian noise, the sample mean already has exponential
concentration and remains a natural and competitive choice \revision{(but it does not outperform the robust alternatives in this setting)}.

\section{Background}
\label{sec:storm-background}
We assume that we have an `ideal' problem that we ultimately wish to solve, 
\begin{align}
    \min_{x\in\mathbb R^d}\phi(x), \label{eq_obj_true}
\end{align}
where $\phi:\R^d\to\R$ is smooth and nonconvex.
However, we also assume that we never have access to $\phi(x)$, and instead we can only calculate
\[
    F(x,\xi) \approx \phi(x),
\]
where \(\xi\) represents the randomness in the evaluation (in some suitable probability space).
If the algorithm requests $N$ evaluations of the objective at the same point, we write
\begin{align}
    Y_i(x)=F(x,\xi_i), \qquad i=1,\ldots,N, \label{eq_noise_generic}
\end{align}
where \(\xi_1,\ldots,\xi_N\) are i.i.d.~copies of \(\xi\).

An \emph{aggregation rule} is any approach to combine these noisy observations into one (hopefully less noisy) estimate:
\begin{align}
    \widehat f(x)
    =
    A_N\bigl(Y_1(x),\ldots,Y_N(x)\bigr). \label{eq_agg_rule_generic}
\end{align}
Here \(A_N:\R^N\to\R\) denotes the finite-sample aggregation rule, such as the sample mean,
$A_N(Y_1,\ldots,Y_N) := \frac{1}{N} \sum_{i=1}^{N} Y_i$.
We assume that, as $N\to\infty$, the aggregation rule converges to some deterministic population target $f(x)$.

Different aggregation rules can define different population targets.
For example, the sample mean aggregation rule gives
\(f_{\mathrm{mean}}(x)=\mathbb E_{\xi}[F(x,\xi)]\),
while using the sample median gives
\(f_{\mathrm{med}}(x)=\mathrm{median}_{\xi}(F(x,\xi))\).
Since we have not defined the specific form of the randomness,
different aggregation rules may yield different population targets.
By taking noisy samples, therefore, we can only really expect
$\widehat{f}(x) \approx f(x)$, rather than $\phi(x)$, and so we can only
expect our DFO algorithm to solve the problem
\begin{align}
    \min_{x\in\R^d} f(x). \label{eq_obj_target}
\end{align}
Thus, we need to assume that solving \eqref{eq_obj_target} would give a
minimizer that is also a solution to our `ideal' problem
\eqref{eq_obj_true}.

\begin{assumption} \label{ass_consistency}
    The noisy objective values $F(x,\xi)$ and aggregation rule $A_N$
    give a population target $f(x)$ such that 
    \[
    f(x_1)\leq f(x_2)
    \quad\Longleftrightarrow\quad
    \phi(x_1)\leq \phi(x_2),
    \qquad
    \forall x_1,x_2\in\mathbb R^d.
    \]
\end{assumption}

For the aggregation rules and noise models studied in this paper, we
identify the relevant population target explicitly and verify when it
satisfies Assumption~\ref{ass_consistency}.
In the main cases considered below, this target is either the population
mean or the population median; for a fixed-proportion trimmed mean,
however, the limiting target may instead be a population trimmed mean.
In our numerical experiments, we compare the aggregation rules by
measuring optimization progress relative to the ideal objective
$\phi(x)$.
Section~\ref{sec:aggregation-bounds} considers specific combinations of
noise models and aggregation rules and discusses their corresponding
population targets and consistency properties.

In the setting where we wish to solve \eqref{eq_obj_target} with access
only to random approximations $\widehat{f}(x)$, a suitable framework is
that of STORM \cite{chen2018stochastic}.
STORM, or Stochastic Trust-Region using Random Models, extends
deterministic trust-region methods to settings where both the local model
and the function estimates used to accept/reject a step are random.
A simplified version of STORM, targeted to the DFO setting, is given in
Algorithm~\ref{alg:storm-simplified}.

\begin{algorithm}[tb]
\caption{DFO STORM algorithm for \eqref{eq_obj_target}}
\label{alg:storm-simplified}
\begin{algorithmic}[1]

\State \textbf{Input:} initial point \(x_0\), initial radius
\(0<\Delta_0\leq\Delta_{\max}\), maximum radius \(\Delta_{\max}>0\),
acceptance parameters \(\eta_1\in(0,1)\), \(\eta_2>0\), and radius
update factor \(\gamma>1\)

\For{\(k=0,1,2,\ldots\)}

    \State Build a quadratic model \(m_k(y) \approx f(y)\) near \(x_k\)
    by interpolation to $p$ independent noisy estimates
    $\widehat{f}(y_1),\ldots,\widehat{f}(y_p)$.

    \State Compute a trial step
    \[
        s_k \approx \argmin_{\|s\|\leq \Delta_k} m_k(x_k+s).
    \]

    \State Get independent noisy estimates $\widehat{f}(x_k)$ and
    $\widehat{f}(x_k+s_k)$, and compute

    \[
        \rho_k
        =
        \frac{
        \widehat f(x_k)-\widehat f(x_k+s_k)
        }{
        m_k(x_k)-m_k(x_k+s_k)
        }.
    \]

    \If{
        \(\rho_k\geq\eta_1\) and
        \(\|\nabla m_k(x_k)\|\geq\eta_2\Delta_k\)
    }
        \State Set \(x_{k+1}=x_k+s_k\) and $\Delta_{k+1}
            =
            \min\{\gamma\Delta_k,\Delta_{\max}\}$.
    \Else
        \State Set \(x_{k+1}=x_k\) and $\Delta_{k+1}
            =
            \gamma^{-1}\Delta_k$.
    \EndIf

\EndFor

\end{algorithmic}
\end{algorithm}

The necessary quality of the model $m_k$ in
Algorithm~\ref{alg:storm-simplified} is given by the following
definition.

\begin{definition}[Fully linear model \cite{conn2009introduction}]
A model \(m_k\) is fully linear for \(f\) on \(B(x_k,\Delta_k)\) if there
exist constants \(\kappa_{ef},\kappa_{eg}>0\) such that
\[
    |f(x_k+s)-m_k(x_k+s)|
    \leq
    \kappa_{ef}\Delta_k^2,
\]
and
\[
    \|\nabla f(x_k+s)-\nabla m_k(x_k+s)\|
    \leq
    \kappa_{eg}\Delta_k,
\]
for all \(\|s\|\leq \Delta_k\).
\end{definition}

In \cite{chen2018stochastic,blanchet2019stochastic}, the convergence and
expected iteration complexity of Algorithm~\ref{alg:storm-simplified}
are given under the following assumptions.

First, we have two standard assumptions on the problem
\eqref{eq_obj_target} and the quality of the step $s_k$.

\begin{assumption}[Smoothness]\label{ass:smooth}
There exists \(f_{\inf}\in\mathbb R\) such that
\[
    f(x)\geq f_{\inf},
    \qquad
    \forall x\in\mathbb R^d.
\]
The function \(f\) is continuously differentiable and has an
\(L\)-Lipschitz continuous gradient for some \(L>0\).
\end{assumption}

\begin{assumption}[Sufficient decrease and bounded model Hessian]
\label{ass:suff-dec}
There exist constants \(\kappa_{scd}>0\) and \(\kappa_H\geq 1\) such that
\(\|H_k\|\leq\kappa_H\) for all \(k\), and the step \(s_k\) satisfies
\[
m_k(x_k)-m_k(x_k+s_k)
\ge
\frac{\kappa_{scd}}{2}\|\nabla m_k(x_k)\|
\min\!\left\{
\frac{\|\nabla m_k(x_k)\|}{\|H_k\|},
\Delta_k
\right\},
\]
where \(H_k=\grad^2 m_k\).
\end{assumption}

Of greater relevance to this work, however, are the following two
assumptions about the accuracy of the model $m_k$ and the estimates
$\widehat{f}(x_k),\widehat{f}(x_k+s_k)$ used to calculate $\rho_k$.

\begin{assumption}[Model accuracy]\label{ass:model-acc}
There exists \(\alpha_m\in(0,1)\) such that, at each iteration \(k\),
\[
    \mathbb P(
        \text{$m_k$ is fully linear for $f$ on $B(x_k,\Delta_k)$}
        \mid \mathcal F_{k-1}
    )
    \geq \alpha_m,
\]
where \(\mathcal F_{k-1}\) is the $\sigma$-algebra generated by all
randomness until the start of iteration $k$, before $m_k$ is constructed.
\end{assumption}

\begin{assumption}[Function estimate accuracy]\label{ass:eval-acc}
There exist \(\alpha_f\in(0,1)\) and \(\kappa_f\geq0\) such that, at each
iteration \(k\),
\[
    \mathbb P\left(
    \max\{
        |\widehat{f}(x_k)-f(x_k)|,\,
        |\widehat{f}(x_k+s_k)-f(x_k+s_k)|
    \}
    \leq
    \kappa_f\Delta_k^2
    \,\middle|\,
    \mathcal F_{k-1/2}
    \right)
    \geq
    \alpha_f,
\]
where \(\mathcal F_{k-1/2}\) is the $\sigma$-algebra generated by all
randomness available after the model \(m_k\) and trial step \(s_k\) have
been constructed, but before the noisy estimates
\(\widehat{f}(x_k)\) and \(\widehat{f}(x_k+s_k)\) used in the acceptance
test are retrieved.
\end{assumption}

The expected iteration complexity of
Algorithm~\ref{alg:storm-simplified} is given by the following result.
To state the explicit bound, let $\eta_1\in(0,1)$ denote the standard
STORM acceptance-ratio threshold and let $\eta_2>0$ denote the threshold
used in the acceptance condition
$\|\nabla m_k(x_k)\|\geq \eta_2\Delta_k$.

\begin{theorem}[Theorem~3.3 of \cite{blanchet2019stochastic}]
Suppose Assumptions~\ref{ass:smooth}--\ref{ass:eval-acc} hold, and
suppose the function-estimate accuracy constant additionally satisfies
\[
    \kappa_f\leq\kappa_{ef},
    \qquad\text{and}\qquad
    \kappa_f
    <
    \frac14\eta_1\eta_2\kappa_{scd}
    \min\left\{\frac{\eta_2}{\kappa_H},1\right\}.
\]
\revision{Suppose also that}
$\alpha_m$ and $\alpha_f$ \revision{are sufficiently close to 1 that} 
\[
    \frac{\alpha_m\alpha_f-\tfrac12}
    {(1-\alpha_m)(1-\alpha_f)}
    \geq
    10+\frac{30L}{40\kappa_{eg}},
\]
and
\[
    \alpha_f
    \geq
    \frac{\kappa_{eg}+0.064L+4\times10^{-4}\eta_2}
    {\kappa_{eg}+0.064L+4.5\times10^{-4}\eta_2}.
\]
Define the first-order stopping time $T_\epsilon
    :=
    \inf\{k\geq0:\|\nabla f(X_k)\|\leq\epsilon\}$.
Then
\[
    \mathbb E[T_\epsilon]
    \leq
    \frac{\alpha_m\alpha_f}{2\alpha_m\alpha_f-1}
    \left(
        \frac{20\Phi_0\kappa_{eg}}{\Theta\epsilon^2}+1
    \right),
\]
where \revision{$\Theta
    =
    \frac{1}{1800}
    \min\{\eta_2\alpha_f,\kappa_{eg}^{-1}\}$ and $\Phi_0
    =
    \nu(f(x_0)-f_{\inf})
    +(1-\nu)\Delta_0^2$,}
\revision{for some $\nu\in(0,1)$.}
\end{theorem}

\revision{In particular, we note that}
larger values of \(\alpha_m\) and \(\alpha_f\) give a smaller
upper bound on the expected number of iterations.

The focus of this work is on the choice of aggregation function, $A_N$.
We shall see that different aggregation functions (and noise models)
require different numbers of samples to achieve sufficiently accurate
approximations, $\widehat{f}(x)\approx f(x)$.
We shall be concerned with determining suitable values of $N$ to achieve
the following requirement.

\begin{definition}
\label{def:prob-accurate}
    A noisy estimate
    \[
        \widehat{f}(x)
        =
        A_N(Y_1(x),\ldots,Y_N(x))
    \]
    is a $(\delta,\beta)$-probabilistically accurate approximation to
    $f(x)$ if
    \begin{align}
        \mathbb P\left(
            |\widehat f(x)-f(x)|\leq\delta
        \right)
        \geq\beta.
        \label{eq_prob_accurate}
    \end{align}
\end{definition}

This notion of accuracy is sufficient to achieve the necessary
probabilistically accurate notions for convergence of STORM.

For the model-construction result below, we consider interpolation models
constructed either by linear interpolation on \(d+1\) suitably poised
points or by \revision{minimum Frobenius norm quadratic interpolation to $p\in\{d+2,\ldots,(d+1)(d+2)/2\}$} 
suitably poised points.

\begin{theorem}
\label{thm_storm_assumptions}
    We have
    \begin{enumerate}[label=(\alph*)]

        \item
        Suppose \(m_k\) is constructed using one of the interpolation
        schemes described above, with \(p\) interpolation points
        \(y_1,\ldots,y_p\).
        If all
        \(\widehat{f}(y_1),\ldots,\widehat{f}(y_p)\)
        are independently
        \((\kappa_f\Delta_k^2,\alpha_m^{1/p})\)-probabilistically
        accurate, then Assumption~\ref{ass:model-acc} holds;
        \label{itm_storm_models}

        \item
        If \(\widehat{f}(x_k)\) and
        \(\widehat{f}(x_k+s_k)\) are independently
        \((\kappa_f\Delta_k^2,\alpha_f^{1/2})\)-probabilistically
        accurate, then Assumption~\ref{ass:eval-acc} holds.

    \end{enumerate}
\end{theorem}

\begin{proof}
Part (a) follows from the interpolation error results in
\cite[Theorems~7.3 and~7.4]{roberts2025interpolation} for linear and
\revision{minimum Frobenius norm}  quadratic interpolation, respectively, provided the
interpolation sets satisfy the corresponding poisedness conditions.
The required joint probability is then obtained from independence of
the \(p\) function estimates.
Part (b) follows directly from independence of the two function
estimates used in the acceptance test.
\end{proof}

In the next section, we assess different choices of aggregation function
$A_N$ in terms of the number of samples required to satisfy
\eqref{eq_prob_accurate}.
In particular, for the finite-variance and failure-noise settings
considered below, robust estimators can achieve a larger value of
\(\beta\) for the same number of samples.
By Theorem~\ref{thm_storm_assumptions}, larger pointwise accuracy
probabilities allow larger values of \(\alpha_m\) and \(\alpha_f\);
and, from the bound on \(T_\epsilon\) above, larger
\(\alpha_m\) and \(\alpha_f\) can reduce the upper bound on the expected
number of iterations.
The model-construction requirement becomes especially demanding when
\(p\) is large, since each of the \(p\) interpolation estimates must be
accurate simultaneously.
Thus, the potential benefit of high-confidence robust estimators is
particularly relevant for models with many interpolation points, such
as those arising in higher-dimensional problems.

\section{Aggregation Bounds for Noisy Function Evaluations}
\label{sec:aggregation-bounds}

We now study the four aggregation rules used in this paper under several
noise models.
For each combination, our main question is how many samples are needed
to satisfy the probabilistic accuracy requirement
\eqref{eq_prob_accurate}.

Fix \(x\in\R^d\), and let \revision{$Y_1(x),\ldots,Y_N(x)$}
be the i.i.d.~noisy evaluations defined in
\eqref{eq_noise_generic}.
We write
\[
    Y_{(1)}(x)\leq\cdots\leq Y_{(N)}(x)
\]
for the corresponding order statistics.

\begin{description}

    \item[Sample Mean]
    The sample mean is
    \[
        \widehat f_{\mathrm{mean}}(x)
        :=
        \frac{1}{N}\sum_{i=1}^N Y_i(x).
    \]
    When the expectation exists, its population target is
    \[
        f_{\mathrm{mean}}(x)
        :=
        \mathbb E_{\xi}[F(x,\xi)].
    \]

    \item[Sample Median]
    Throughout this paper, we use the usual midpoint convention for the
    sample median:
    \[
        \widehat f_{\mathrm{med}}(x)
        :=
        \begin{cases}
            Y_{((N+1)/2)}(x),
            & \text{if \(N\) is odd}, \\[1mm]
            \dfrac{
                Y_{(N/2)}(x)+Y_{(N/2+1)}(x)
            }{2},
            & \text{if \(N\) is even}.
        \end{cases}
    \]
    When \(F(x,\xi)\) has a unique population median, we denote this
    population target by
    \[
        f_{\mathrm{med}}(x)
        :=
        \operatorname{median}_{\xi}(F(x,\xi)).
    \]

    \item[Trimmed Mean]
    For an integer
    \[
        r\in\{0,\ldots,\lfloor(N-1)/2\rfloor\},
    \]
    the two-sided \(r\)-trimmed mean is
    \[
        \widehat f_{\mathrm{trim},r}(x)
        :=
        \frac{1}{N-2r}
        \sum_{i=r+1}^{N-r}Y_{(i)}(x).
    \]
    The choice of \(r\) depends on the setting.
    For the finite-variance analysis, \(r\) is chosen according to the
    required confidence level.
    For failure noise, we consider a fixed trimming proportion.
    The numerical experiments use
    \(r=\lfloor0.1N\rfloor\).
    The population target of a trimmed mean can depend on how \(r\)
    scales with \(N\), as discussed where this distinction is needed
    below.

    \item[Median-of-Means]
    For the theoretical results below, let \(N=Km\).
    We divide the observations into \(K\) blocks
    \(B_1,\ldots,B_K\), each containing \(m\) observations.
    The corresponding block means are
    \[
        \bar Y_j(x)
        :=
        \frac{1}{m}\sum_{i\in B_j}Y_i(x),
        \qquad
        j=1,\ldots,K.
    \]
    The median-of-means estimator is
    \[
        \widehat f_{\mathrm{MoM}}(x)
        :=
        \operatorname{median}
        \{\bar Y_1(x),\ldots,\bar Y_K(x)\}.
    \]
    In the settings considered below, it is used to estimate the mean
    population target \(f_{\mathrm{mean}}(x)\).
    The theoretical choice of \(K\) depends on the required confidence
    level.
    In the numerical experiments, we instead use \(K=5\).
    When \(N\) is not divisible by \(5\), the implementation uses
    nearly equal-sized blocks, as described in
    Section~\ref{sec:numerical}.

\end{description}

The sample median, trimmed mean, and median-of-means reduce the influence
of extreme observations relative to the sample mean.
Their guarantees depend on the noise model.
We consider three settings.

For the first two settings, we use the additive-noise model
\[
    Y(x)=\phi(x)+\xi,
    \qquad
    \mathbb E[\xi]=0.
\]

\begin{assumption}[Sub-Gaussian noise]
\label{ass_sub_gaussian}
The noise \(\xi\) is sub-Gaussian with variance proxy \(\sigma^2\), that is,
\[
    \mathbb E[e^{t\xi}]
    \leq
    e^{\sigma^2t^2/2},
    \qquad
    \forall t\in\R.
\]
\end{assumption}

\begin{assumption}[Finite variance noise]
\label{ass_finite_variance}
The noise \(\xi\) satisfies \revision{$\mathrm{Var}(\xi)\leq\sigma^2<\infty$.}
\end{assumption}

A sub-Gaussian random variable with variance proxy \(\sigma^2\), as
defined in Assumption~\ref{ass_sub_gaussian}, satisfies \liurevision{$\mathrm{Var}(\xi)\leq\sigma^2$ by \cite[Eq.~(3.2)]{zhang2021concentration}}.
Thus, Assumption~\ref{ass_sub_gaussian} is a special case of
Assumption~\ref{ass_finite_variance}.

Under either Assumption~\ref{ass_sub_gaussian} or
Assumption~\ref{ass_finite_variance}, the mean population target is
\[
    f_{\mathrm{mean}}(x)=\phi(x).
\]
If \(\xi\) also has a unique population median \(m_\xi\), then
\[
    f_{\mathrm{med}}(x)
    =
    \phi(x)+m_\xi.
\]
Since \(m_\xi\) does not depend on \(x\), this median target also satisfies
Assumption~\ref{ass_consistency}.

We also consider stochastic failures in the objective evaluation.

\begin{assumption}[Failure noise]
\label{ass_failure}
The noisy evaluation satisfies
\[
Y(x)=
\begin{cases}
    \phi(x), & \text{with probability }q,\\
    V, & \text{with probability }1-q,
\end{cases}
\]
for some \(q\in(1/2,1]\) and \(V\in\R\) independent of \(x\).
\end{assumption}

Under Assumption~\ref{ass_failure},
\[
    f_{\mathrm{mean}}(x)
    =
    q\phi(x)+(1-q)V,
    \qquad
    f_{\mathrm{med}}(x)
    =
    \phi(x).
\]
Both targets satisfy Assumption~\ref{ass_consistency}.
For later use, define
\[
    D(x):=|V-\phi(x)|.
\]

\subsection{Sample mean}
\label{subsec:sample-mean}

We first consider the sample mean under finite-variance noise.
We then note the sharper result available under the stronger
sub-Gaussian assumption.

\subsubsection{Finite-variance and sub-Gaussian noise}

\revision{We now give sample estimates required for the sample mean to be sufficiently accurate.}

\begin{proposition}[Sample mean under finite-variance noise]
\label{prop:mean-chebyshev}
Under Assumption~\ref{ass_finite_variance},
for any
\(\delta>0\) and \(\beta\in(0,1)\), the sample mean is a
\((\delta,\beta)\)-probabilistically accurate approximation to
\(f_{\mathrm{mean}}(x)\) provided that
\[
    N
    \geq
    \frac{\sigma^2}{\delta^2(1-\beta)}.
\]
\end{proposition}
\begin{proof}
    \revision{This is an immediate consequence of Chebyshev's inequality, \cite[Corollary~1.2.5]{vershynin2018high}.}
\end{proof}

Under only a finite-variance assumption, the sufficient sample size
therefore has \(O((1-\beta)^{-1})\) dependence on the inverse failure
probability.

A stronger guarantee is available under sub-Gaussian noise.

\begin{proposition}[Sample mean under sub-Gaussian noise]
\label{prop:mean-subgaussian}
Under Assumption~\ref{ass_sub_gaussian}, 
for any
\(\delta>0\) and \(\beta\in(0,1)\), the sample mean is a
\((\delta,\beta)\)-probabilistically accurate approximation to
\(f_{\mathrm{mean}}(x)\) provided that
\[
    N
    \geq
    \frac{2\sigma^2}{\delta^2}
    \log\!\left(\frac{2}{1-\beta}\right).
\]
\end{proposition}
\begin{proof}
    \revision{This is an immediate consequence of Hoeffding's inequality, \cite[Proposition~2.5]{wainwright2019high}.}
\end{proof}

Thus, the dependence on the inverse failure probability improves from
\(O((1-\beta)^{-1})\) under only finite variance to
\(O(\log(1/(1-\beta)))\) under sub-Gaussian noise.
This difference becomes important when a very high confidence level is
required.

\subsubsection{Failure noise}

\begin{proposition}[Sample mean under failure noise]
\label{prop:mean-failure}
Under Assumption~\ref{ass_failure}, suppose \(D(x)>0\).
\revision{Then}
the sample mean is a
\((\delta,\beta)\)-probabilistically accurate approximation to
\(f_{\mathrm{mean}}(x)\) provided that
\[
    N
    \geq
    \frac{D(x)^2}{2\delta^2}
    \log\!\left(\frac{2}{1-\beta}\right).
\]
If \(D(x)=0\), every evaluation is equal to \(\phi(x)\), so the result
holds trivially for every \(N\geq1\).
\end{proposition}
\begin{proof}
    \revision{This follows from Hoeffding's inequality for bounded random variables, \cite[Theorem~2.2.6]{vershynin2018high}.}
\end{proof}

The required sample size depends on \(D(x)\).
Thus, estimating the mean target becomes more expensive when the failed
value is far from the noise-free objective value.

The mean population target itself satisfies
\[
    |f_{\mathrm{mean}}(x)-\phi(x)|
    =
    (1-q)D(x).
\]
This difference does not prevent the optimization of \(\phi\) under
Assumption~\ref{ass_consistency}.
It does show, however, that the function estimates need not be close to
the noise-free values when \(D(x)\) is large.

\subsection{Trimmed mean}
\label{subsec:trimmed-mean}

We next derive finite-variance and failure-noise guarantees for the
trimmed mean.

\subsubsection{Finite-variance noise}

The following finite-sample bound follows from
\cite[Theorem~1.1.1]{oliveira2025finite}.

\begin{theorem}[Trimmed mean under finite variance]
\label{thm:tm-finite-variance}
Suppose Assumption~\ref{ass_finite_variance} holds.
Let
\[
    c_0=3\sqrt2+8,
    \qquad
    c_1=4+4\sqrt2.
\]
Let \(t>0\) satisfy
\[
    0<t\leq
    \sqrt{
        \frac{N}{(\sqrt2+1)^2}-2
    },
\]
and choose
\[
    r(t)
    =
    \left\lceil\frac{t^2}{2}\right\rceil.
\]
Then
\[
    \mathbb P\left(
        \left|
        \widehat f_{\mathrm{trim},r(t)}(x)
        -
        f_{\mathrm{mean}}(x)
        \right|
        >
        \frac{(c_0+c_1t)\sigma}{\sqrt N}
    \right)
    \leq
    4\exp\!\left(-\frac{t^2}{2}\right).
\]
\end{theorem}

\begin{corollary}[Probabilistic accuracy of the trimmed mean]
\label{cor:tm-finite-variance}
Under Assumption~\ref{ass_finite_variance}, let
\[
    t_\beta
    :=
    \sqrt{
        2\log\!\left(\frac{4}{1-\beta}\right)
    },
\]
and choose
\[
    r_\beta
    :=
    \left\lceil\frac{t_\beta^2}{2}\right\rceil
    =
    \left\lceil
        \log\!\left(\frac{4}{1-\beta}\right)
    \right\rceil.
\]
If
\[
    N
    \geq
    \frac{(c_0+c_1t_\beta)^2\sigma^2}{\delta^2},
\]
and the range condition on \(t_\beta\) in
Theorem~\ref{thm:tm-finite-variance} holds, then
\(\widehat f_{\mathrm{trim},r_\beta}(x)\) is a
\((\delta,\beta)\)-probabilistically accurate approximation to
\(f_{\mathrm{mean}}(x)\).
\end{corollary}

For fixed \(\beta\), the value \(r_\beta\) does not grow with \(N\).
Hence,
\[
    \frac{r_\beta}{N}\to0
    \qquad\text{as }N\to\infty.
\]
Up to constants, Corollary~\ref{cor:tm-finite-variance} gives
\[
    N
    =
    O\left(
        \frac{\sigma^2}{\delta^2}
        \log\frac{1}{1-\beta}
    \right).
\]
Thus, under only finite variance, the trimmed mean has logarithmic
dependence on the inverse failure probability.
The Chebyshev guarantee for the sample mean instead has
\(1/(1-\beta)\) dependence.

\subsubsection{Failure noise}

We now consider the trimmed mean under Assumption~\ref{ass_failure}.

\begin{proposition}[Exact recovery under failure noise]
\label{prop:tm-exact-recovery}
Under Assumption~\ref{ass_failure}, let \(r\) be an admissible trimming
level as defined above.
If
\[
    \frac{r}{N}>1-q,
\]
then
\[
    \mathbb P\left(
        \widehat f_{\mathrm{trim},r}(x)=\phi(x)
    \right)
    \geq
    1-
    \exp\!\left(
        -2N
        \left(
            \frac{r}{N}-(1-q)
        \right)^2
    \right).
\]
\end{proposition}

\begin{proof}
Let
\[
    M_N
    :=
    \sum_{i=1}^{N}
    \mathbf 1\{Y_i(x)=V\}
\]
denote the number of failed evaluations.
Then
\[
    M_N\sim\mathrm{Bin}(N,1-q).
\]

If \(V=\phi(x)\), the result is immediate.
Suppose instead that \(V\neq\phi(x)\).
All failed observations then lie at one end of the ordered sample.
On the event \(\{M_N\leq r\}\), the trimming step removes every failed
observation.
Therefore,
\[
    \widehat f_{\mathrm{trim},r}(x)=\phi(x).
\]
It follows that
\[
    \left\{
        \widehat f_{\mathrm{trim},r}(x)\neq\phi(x)
    \right\}
    \subseteq
    \{M_N>r\}.
\]

Since \(M_N\sim\mathrm{Bin}(N,1-q)\), Hoeffding's inequality
\cite[Theorem~2.2.6]{vershynin2018high} gives
\[
    \mathbb P(M_N>r)
    \leq
    \exp\!\left(
        -2N
        \left(
            \frac{r}{N}-(1-q)
        \right)^2
    \right),
\]
whenever \(r/N>1-q\).
\end{proof}

If \(r=\lfloor\tau N\rfloor\), the condition becomes
\[
    \frac{\lfloor\tau N\rfloor}{N}>1-q.
\]
For large \(N\), this is essentially
\[
    \tau>1-q.
\]
Thus, the trimming proportion must exceed the failure probability for
the exact-recovery result above to apply.

This condition also shows the difference between the theorem and our
numerical implementation.
The numerical experiments use \(\tau=0.1\), while the failure
probability is \(1-q=0.2\).
Therefore, Proposition~\ref{prop:tm-exact-recovery} does not predict
exact recovery for the trimmed mean used in those experiments.
\liurevision{In the numerical experiments, the failed evaluations all return the
fixed value \(V=-10^{4}\), as in Assumption~\ref{ass_failure}.}

\revision{More general contamination models have also been studied.
For example, \cite[Theorem~1.3.1]{oliveira2025finite} gives, under suitable assumptions, a concentration inequality for a failure noise model where the `failure' values can be arbitrary (not all equal as required by Assumption~\ref{ass_failure}), and the `successful' values can also be noisy.}

\subsection{Median-of-means}
\label{subsec:mom}

We next consider the median-of-means estimator.
For the finite-sample theory below, \(K\) is chosen according to the
required confidence level and is taken to be odd.

\subsubsection{Finite-variance noise}
\label{subsec:mom-finite-variance}

\begin{theorem}[Median-of-means under finite variance]
\label{thm:mom-finite-variance}
Suppose Assumption~\ref{ass_finite_variance} holds.
Let \(N=Km\), where \(K\) is odd.
\revision{Then}
\[
    \mathbb P\left(
        \left|
        \widehat f_{\mathrm{MoM}}(x)-f_{\mathrm{mean}}(x)
        \right|
        >
        2\sigma\sqrt{\frac{K}{N}}
    \right)
    \leq
    \exp\left(-\frac{K}{8}\right).
\]
\end{theorem}
\begin{proof}
    \revision{See \cite[Theorem~2]{lugosi2019mean}.}
\end{proof}

\begin{corollary}[Probabilistic accuracy of median-of-means]
\label{cor:mom-finite-variance}
Under Assumption~\ref{ass_finite_variance}, choose an odd integer \(K\)
such that
\[
    K
    \geq
    8\log\left(\frac{1}{1-\beta}\right),
\]
and let \(N=Km\) satisfy
\[
    N
    \geq
    \frac{4\sigma^2K}{\delta^2}.
\]
Then
\(\widehat f_{\mathrm{MoM}}(x)\) is a
\((\delta,\beta)\)-probabilistically accurate approximation to
\(f_{\mathrm{mean}}(x)\).
\end{corollary}

Up to constants, this gives
\[
    N
    =
    O\left(
        \frac{\sigma^2}{\delta^2}
        \log\frac{1}{1-\beta}
    \right).
\]
Thus, median-of-means has logarithmic dependence on the inverse failure
probability under only finite variance.
This improves on the \(1/(1-\beta)\) dependence of the Chebyshev
guarantee in Proposition~\ref{prop:mean-chebyshev}.

Since sub-Gaussian noise is a special case of finite-variance noise,
Theorem~\ref{thm:mom-finite-variance} also applies under
Assumption~\ref{ass_sub_gaussian}.
In that stronger setting, however, it does not improve the sample-size
order already achieved by the sample mean.

\subsubsection{Failure noise}
\label{subsec:mom-failure}

Under Assumption~\ref{ass_failure}, median-of-means targets the same
mean population target as the sample mean.
The main question is therefore how many samples are needed to estimate
this target with high probability.

\begin{proposition}[Median-of-means under failure noise]
\label{prop:mom-failure}
Let \(N=Km\), with \(K\) odd, and suppose \(D(x)>0\).
Define
\[
    p_m(\delta)
    :=
    2\exp\left(
        -\frac{2m\delta^2}{D(x)^2}
    \right).
\]
If
\[
    p_m(\delta)<\frac12,
\]
then
\[
    \mathbb P\left(
        \left|
        \widehat f_{\mathrm{MoM}}(x)-f_{\mathrm{mean}}(x)
        \right|
        >
        \delta
    \right)
    \leq
    \exp\left(
        -2K
        \left(
            \frac12-p_m(\delta)
        \right)^2
    \right).
\]
If \(D(x)=0\), the estimator is exact.
\end{proposition}

\begin{proof}
For each block \(B_j\), the observations lie in the interval with
endpoints \(\phi(x)\) and \(V\).
This interval has width \(D(x)\).
Hoeffding's inequality
\cite[Theorem~2.2.6]{vershynin2018high} therefore gives
\[
    \mathbb P\left(
        |\bar Y_j(x)-f_{\mathrm{mean}}(x)|>\delta
    \right)
    \leq
    p_m(\delta).
\]

Define
\[
    Z_j
    :=
    \mathbf 1\{
        |\bar Y_j(x)-f_{\mathrm{mean}}(x)|>\delta
    \}.
\]
The variables \(Z_1,\ldots,Z_K\) are independent and satisfy
\[
    \mathbb P(Z_j=1)\leq p_m(\delta).
\]
If
\[
    |\widehat f_{\mathrm{MoM}}(x)-f_{\mathrm{mean}}(x)|>\delta,
\]
then at least half of the block means must be bad.
Hence,
\[
    \left\{
        |\widehat f_{\mathrm{MoM}}(x)-f_{\mathrm{mean}}(x)|>\delta
    \right\}
    \subseteq
    \left\{
        \sum_{j=1}^K Z_j\geq\frac K2
    \right\}.
\]
Since \(p_m(\delta)<1/2\), another application of Hoeffding's inequality
gives
\[
    \mathbb P\left(
        \sum_{j=1}^K Z_j\geq\frac K2
    \right)
    \leq
    \exp\left(
        -2K
        \left(
            \frac12-p_m(\delta)
        \right)^2
    \right).
\]
\end{proof}

A simple sufficient choice is obtained by requiring
\[
    p_m(\delta)\leq\frac14.
\]
It is enough to take
\[
    m
    \geq
    \frac{D(x)^2}{2\delta^2}\log 8
\]
and to choose an odd integer \(K\) satisfying
\[
    K
    \geq
    8\log\left(\frac{1}{1-\beta}\right).
\]
Therefore,
\[
    N=Km
    =
    O\left(
        \frac{D(x)^2}{\delta^2}
        \log\frac{1}{1-\beta}
    \right).
\]

Thus, under the fixed-value failure model, median-of-means can estimate
its population target with high probability.
Its required block size, however, depends quadratically on \(D(x)\).
The fact that its population target differs from the noise-free value
\(\phi(x)\) is not itself an obstacle when
Assumption~\ref{ass_consistency} holds.
The more important point is that the required sample size depends on
both the desired accuracy \(\delta\) and the magnitude of the failed
value.

\subsection{Sample median}
\label{subsec:sample-median}

We finally consider the sample median.
Unlike the robust mean estimators above, its natural population target
is \(f_{\mathrm{med}}(x)\).
Its concentration depends on how much probability mass lies near this
median.

\begin{assumption}[Median local mass]
\label{ass:median-target}
For the fixed \(x\) above and a given accuracy level \(\delta>0\), let
\(G_x\) denote the cumulative distribution function of \(Y(x)\),
and suppose that \(Y(x)\) has a unique population median
\(f_{\mathrm{med}}(x)\).

Define
\[
    \gamma_x(\delta)
    :=
    \min\left\{
        G_x(f_{\mathrm{med}}(x)+\delta)-\frac12,\,
        \frac12-G_x(f_{\mathrm{med}}(x)-\delta)
    \right\}.
\]
We assume that
\[
    \gamma_x(\delta)>0.
\]
\end{assumption}

The quantity \(\gamma_x(\delta)\) measures how far the CDF has moved
away from \(1/2\) at distance \(\delta\) on either side of the median.
It is therefore a local condition on the distribution near its median.
It does not require \(G_x\) to be continuous.
It also does not require a finite mean or finite variance.

The Dvoretzky--Kiefer--Wolfowitz \revision{(DKW)} inequality used below is valid for
arbitrary distribution functions, including discontinuous ones; see
\cite[Theorem~11.5]{kosorok2008introduction}.
See also \cite[Lemma~7.1]{morvai2024countable}, which discusses the
extension from continuous to possibly discontinuous distribution
functions.
The sharp constant in the DKW inequality is due to
Massart~\cite{massart1990tight}.

\subsubsection{General concentration via DKW}

Let \(G_{x,N}\) denote the empirical cumulative distribution function
of \(Y_1(x),\ldots,Y_N(x)\).
The midpoint estimator defined above is a sample median of this empirical
distribution.
Hence,
\[
    G_{x,N}(\widehat f_{\mathrm{med}}(x)^-)
    \leq
    \frac12
    \leq
    G_{x,N}(\widehat f_{\mathrm{med}}(x)),
\]
\revision{where we use the notation $f(x^{-}) := \lim_{y\uparrow x} f(y)$.}

\begin{proposition}[Median concentration via DKW]
\label{prop:median-dkw}
Suppose Assumption~\ref{ass:median-target} holds at \(x\) for the
accuracy level \(\delta>0\).
Then
\[
    \mathbb P\left(
        |\widehat f_{\mathrm{med}}(x)-f_{\mathrm{med}}(x)|>\delta
    \right)
    \leq
    2\exp\!\left(
        -2N\gamma_x(\delta)^2
    \right).
\]
\end{proposition}

\begin{proof}
First suppose that \revision{$\widehat f_{\mathrm{med}}(x)
    >
    f_{\mathrm{med}}(x)+\delta$}.
Since \(G_{x,N}\) is nondecreasing,
\[
    G_{x,N}(f_{\mathrm{med}}(x)+\delta)
    \leq
    G_{x,N}(\widehat f_{\mathrm{med}}(x)^-)
    \leq
    \frac12.
\]
By the definition of \(\gamma_x(\delta)\),
\[
    G_x(f_{\mathrm{med}}(x)+\delta)-\frac12
    \geq
    \gamma_x(\delta).
\]
Therefore,
\[
    G_x(f_{\mathrm{med}}(x)+\delta)
    -
    G_{x,N}(f_{\mathrm{med}}(x)+\delta)
    \geq
    \gamma_x(\delta).
\]

Now suppose that \revision{$\widehat f_{\mathrm{med}}(x)
    <
    f_{\mathrm{med}}(x)-\delta$}.
Then
\[
    G_{x,N}(f_{\mathrm{med}}(x)-\delta)
    \geq
    G_{x,N}(\widehat f_{\mathrm{med}}(x))
    \geq
    \frac12.
\]
Also,
\[
    \frac12-G_x(f_{\mathrm{med}}(x)-\delta)
    \geq
    \gamma_x(\delta).
\]
Hence,
\[
    G_{x,N}(f_{\mathrm{med}}(x)-\delta)
    -
    G_x(f_{\mathrm{med}}(x)-\delta)
    \geq
    \gamma_x(\delta).
\]

Combining the two cases gives
\[
    \left\{
        |\widehat f_{\mathrm{med}}(x)-f_{\mathrm{med}}(x)|>\delta
    \right\}
    \subseteq
    \left\{
        \sup_{t\in\mathbb R}
        |G_{x,N}(t)-G_x(t)|
        \geq
        \gamma_x(\delta)
    \right\}.
\]

For any \(0<\varepsilon<\gamma_x(\delta)\),
\[
    \left\{
        \sup_{t\in\mathbb R}|G_{x,N}(t)-G_x(t)|
        \geq\gamma_x(\delta)
    \right\}
    \subseteq
    \left\{
        \sup_{t\in\mathbb R}|G_{x,N}(t)-G_x(t)|
        >\varepsilon
    \right\}.
\]
The DKW inequality
\cite[Theorem~11.5]{kosorok2008introduction} therefore gives
\[
    \mathbb P\left(
        \sup_{t\in\mathbb R}|G_{x,N}(t)-G_x(t)|
        \geq\gamma_x(\delta)
    \right)
    \leq
    2\exp(-2N\varepsilon^2).
\]
Letting \(\varepsilon\uparrow\gamma_x(\delta)\) gives the stated bound.
\end{proof}

Consequently, the sample median is a
\((\delta,\beta)\)-probabilistically accurate approximation to
\(f_{\mathrm{med}}(x)\) if
\[
    N
    \geq
    \frac{1}{2\gamma_x(\delta)^2}
    \log\!\left(\frac{2}{1-\beta}\right).
\]
The dependence on the confidence level is therefore logarithmic in
\(1/(1-\beta)\).
The remaining quantity, \(\gamma_x(\delta)\), depends on the local
shape of the distribution.
If little probability mass lies near the median, then
\(\gamma_x(\delta)\) is small and more samples are needed.

\paragraph{A sufficient local-density condition.}

The dependence of \(\gamma_x(\delta)\) on \(\delta\) can be made more
explicit under a local density assumption.
Suppose that \(G_x\) has a density \(g_x\) in a neighbourhood of
\(f_{\mathrm{med}}(x)\).
Suppose further that there exist constants
\(\underline g_x>0\) and \(\delta_0>0\) such that
\[
    g_x(t)\geq\underline g_x
    \qquad
    \text{whenever }
    |t-f_{\mathrm{med}}(x)|\leq\delta_0.
\]
Then, for every \(0<\delta\leq\delta_0\),
\[
    \gamma_x(\delta)
    \geq
    \underline g_x\,\delta.
\]
It follows that
\[
    N
    \geq
    \frac{1}{2\underline g_x^2\delta^2}
    \log\!\left(\frac{2}{1-\beta}\right)
\]
is sufficient.

For a location-scale family with scale \(\sigma\), suppose that the
standardized density is bounded below by a positive constant in a
neighbourhood of its median.
\revision{That is,} \(\underline g_x\) is proportional to \(\sigma^{-1}\).
The sample-size requirement therefore has the familiar order
\[
    N
    =
    O\left(
        \frac{\sigma^2}{\delta^2}
        \log\frac{1}{1-\beta}
    \right).
\]
Gaussian noise is one example for which such a local density condition
holds.

Finite variance alone does not give a quantitative lower bound on the
probability mass near the median.
In general, the concentration of the sample median is therefore
controlled by \(\gamma_x(\delta)\), rather than by the variance alone.

\subsubsection{Failure noise}

Under Assumption~\ref{ass_failure}, a successful evaluation occurs with
probability \(q>1/2\).
A strict majority of successful observations is therefore enough for
the sample median to recover the noise-free value exactly.

The following result also covers a slightly more general setting.
The failed evaluations do not need to take the same value \(V\), \revision{as required by Assumption~\ref{ass_failure}}.

\begin{theorem}[Exact recovery under majority correctness]
\label{thm:median-failure}
Suppose each evaluation returns \(\phi(x)\) independently with
probability \(q>1/2\).
The remaining evaluations may take arbitrary values.
Then
\[
    \mathbb P\left(
        \widehat f_{\mathrm{med}}(x)=\phi(x)
    \right)
    \geq
    1-
    \exp\!\left(
        -2N\left(q-\frac12\right)^2
    \right).
\]
Consequently, for every \(\delta>0\),
\[
    \mathbb P\left(
        |\widehat f_{\mathrm{med}}(x)-f_{\mathrm{med}}(x)|
        \leq\delta
    \right)
    \geq
    1-
    \exp\!\left(
        -2N\left(q-\frac12\right)^2
    \right).
\]
In particular, the sample median is
\((\delta,\beta)\)-probabilistically accurate for every \(\delta>0\)
whenever
\[
    N
    \geq
    \frac{1}{2(q-\frac12)^2}
    \log\!\left(\frac{1}{1-\beta}\right).
\]
\end{theorem}

\begin{proof}
Let
\[
    Z_i
    :=
    \mathbf 1\{\text{the \(i\)th evaluation is successful}\}.
\]
Then \(Z_1,\ldots,Z_N\) are independent Bernoulli random variables with \revision{$\mathbb E[Z_i]=q$}.
If \revision{$\sum_{i=1}^N Z_i>\frac N2$,}
then more than half of the observations equal \(\phi(x)\) \revision{and hence $\widehat f_{\mathrm{med}}(x)=\phi(x)$}.
Therefore,
\[
    \left\{
        \widehat f_{\mathrm{med}}(x)\neq\phi(x)
    \right\}
    \subseteq
    \left\{
        \frac1N\sum_{i=1}^N Z_i\leq\frac12
    \right\}.
\]
Since \(q>1/2\), Hoeffding's inequality
\cite[Theorem~2.2.6]{vershynin2018high} gives
\[
    \mathbb P\left(
        \frac1N\sum_{i=1}^N Z_i\leq\frac12
    \right)
    \leq
    \exp\!\left(
        -2N\left(q-\frac12\right)^2
    \right).
\]
Exact recovery gives zero estimation error.
It therefore implies the stated
\((\delta,\beta)\)-accuracy guarantee for every \(\delta>0\).
\end{proof}

This behaviour is quite different from mean-based estimation under
failure noise.
The median does not average the failed values.
Its sample-size guarantee therefore does not depend on the desired
accuracy \(\delta\) or on the magnitude of the corrupted values.
It only requires successful evaluations to form a majority.

\subsection{Summary of aggregation guarantees}
\label{subsec:aggregation-summary}

The main difference between the estimators under finite-variance noise
is their dependence on the required confidence level.
For the sample mean, Proposition~\ref{prop:mean-chebyshev} gives
\[
    N
    =
    O\left(
        \frac{\sigma^2}{\delta^2}
        \frac{1}{1-\beta}
    \right).
\]
The trimmed mean and median-of-means instead give
\[
    N
    =
    O\left(
        \frac{\sigma^2}{\delta^2}
        \log\frac{1}{1-\beta}
    \right).
\]
Thus, all three mean estimators have the same
\(O(\sigma^2/\delta^2)\) accuracy factor.
The difference lies in the confidence term.
Under only finite variance, the robust estimators replace
\(1/(1-\beta)\) by \(\log(1/(1-\beta))\).

If the noise is sub-Gaussian, the sample mean already has logarithmic
confidence dependence.
The robust mean estimators therefore do not improve this sample-size
order in the sub-Gaussian setting.

The sample median also has logarithmic confidence dependence:
\[
    N
    =
    O\left(
        \frac{1}{\gamma_x(\delta)^2}
        \log\frac{1}{1-\beta}
    \right).
\]
Its accuracy dependence is controlled by the probability mass near the
population median rather than by the variance.
Under the local-density condition above, with density of order
\(\sigma^{-1}\) near the median, this becomes
\[
    N
    =
    O\left(
        \frac{\sigma^2}{\delta^2}
        \log\frac{1}{1-\beta}
    \right).
\]

The confidence term can matter when a model is constructed from many
independently estimated function values.
As discussed in Theorem~\ref{thm_storm_assumptions}, each pointwise
estimate must then be accurate with a probability that can be very close
to one.

\revision{The impact of the confidence term can be material.}
For example, a fully determined quadratic interpolation model in
dimension \(d=10\) uses \revision{$p=66$}
interpolation points.
If the desired model-accuracy probability is \revision{relatively mild, say}
\(\alpha_m=0.9\), then \revision{for Assumption~\ref{ass:model-acc} to hold, Theorem~\ref{thm_storm_assumptions}(a) requires $\beta=\alpha_m^{1/p}\approx0.998405$ for each interpolation point.}
In this case,
\[
    \frac{1}{1-\beta}
    \approx
    6.27\times10^2,
    \qquad \revision{\text{and}} \qquad
    \log\!\left(\frac{1}{1-\beta}\right)
    \approx
    6.44.
\]
The confidence factor in the Chebyshev sample-mean guarantee is therefore
almost two orders of magnitude larger than the logarithmic factor in the
robust bounds.
This example illustrates why the difference can become important when
many accurate interpolation values are needed at the same iteration.

Failure noise leads to a different comparison.
Under the fixed-value failure model, the sample mean requires
\[
    N
    =
    O\left(
        \frac{D(x)^2}{\delta^2}
        \log\frac{1}{1-\beta}
    \right)
\]
samples.
Median-of-means has the same dependence on \(D(x)\), \(\delta\), and
the confidence level, up to constants.

The sample median behaves differently.
It gives exact recovery with
\[
    N
    =
    O\left(
        \log\frac{1}{1-\beta}
    \right),
\]
independently of both \(D(x)\) and \(\delta\).

The trimmed mean has a similar exact-recovery property when the trimming
proportion is large enough.
If a fixed proportion \(\tau\) satisfies
\[
    \tau>1-q,
\]
then Proposition~\ref{prop:tm-exact-recovery} gives exact recovery with
\[
    N
    =
    O\left(
        \frac{1}
        {(\tau-(1-q))^2}
        \log\frac{1}{1-\beta}
    \right).
\]
This bound is also independent of the magnitude of the failed value and
the requested function-value accuracy.

These results highlight two different reasons for using robust
aggregation in stochastic model-based DFO.
Under finite-variance noise, robust mean estimators improve the
high-confidence sample-size dependence relative to the Chebyshev
guarantee for the sample mean.
Under failure noise, the median and a suitably trimmed mean can instead
use the structure of the failures to recover the noise-free objective
value exactly with high probability.

\section{Numerical Experiments}
\label{sec:numerical}

We now investigate how the choice of aggregation rule affects the
performance of a model-based DFO solver.
We consider two test sets.
The first is the low-dimensional Mor\'e--Wild collection.
The second consists of \liurevision{medium-dimensional} problems from CUTEst.
This allows us to examine whether the behaviour observed on the
Mor\'e--Wild problems persists on a more challenging test set.

We compare the sample mean, sample median, trimmed mean, and
median-of-means under four noise models.
We also compare fixed sampling with a radius-dependent sampling rule.

\subsection{Experimental setup}

\paragraph{Test problems and solver.}

The low-dimensional test set consists of the 53 problems from the
Mor\'e--Wild collection \cite{more2009benchmarking}.

For the larger-scale experiments, we use 29 problem configurations from
the CUTEst collection \liurevision{\cite{gould2015cutest}}.
\liurevision{The selection is based on the general-objective
test set in \cite[Appendix~D]{cartis2019improving}.}
The CUTEst problems are instantiated through PyCUTEst
\liurevision{\cite{fowkes2022pycutest}}.
\liurevision{The CUTEst configurations have dimensions
$64\leq d\leq110$, with most having $d=100$,
compared with $2\leq d\leq12$ for the Mor\'e--Wild problems.}
We therefore refer to this collection as the
\liurevision{medium-dimensional} test set below.

Both test sets are solved using \textsc{Py-BOBYQA}
\cite{cartis2019improving}.
The solver is modified so that repeated evaluations at the same point
can be combined using different aggregation rules.

For the Mor\'e--Wild experiments, we set the solver to treat the
objective as noisy and otherwise use the default Py-BOBYQA settings for
noisy problems.
For the CUTEst experiments, we use
\[
    n_{\mathrm{pt}}=2d+1
\]
interpolation points, where \(d\) is the problem dimension.
Bounds supplied by CUTEst are passed to Py-BOBYQA when available.
Within each test set, all aggregation rules use the same solver
settings.



\paragraph{Aggregation rules.}

We compare the four aggregation rules defined in
Section~\ref{sec:aggregation-bounds}: the sample mean, sample median,
trimmed mean, and median-of-means.

For the trimmed mean, we use
\[
    r=\lfloor0.1N\rfloor.
\]
Thus, up to rounding, the smallest and largest \(10\%\) of the
observations are removed.

For median-of-means, we use \(K=5\) blocks.
When \(N\) is not divisible by \(5\), the final block contains the
remaining observations.
If \(N<5\), the implementation uses the sample median as a fallback.

These are fixed practical choices for the numerical experiments.
They differ from some of the confidence-dependent choices used in the
theoretical analysis in Section~\ref{sec:aggregation-bounds}.

\paragraph{Noise models.}

We use the same four noise settings on both test sets.
Throughout this section, \(\phi(x)\) denotes the noise-free objective
value.

\begin{itemize}

    \item \textbf{Gaussian noise.}
    We use
    \[
        F(x,\xi)
        =
        \phi(x)+\xi,
        \qquad
        \xi\sim\mathcal N(0,10^{-4}).
    \]
    The noise therefore has standard deviation \(10^{-2}\).
    This is a sub-Gaussian instance of
    Assumption~\ref{ass_sub_gaussian}.

    \item \textbf{Student-\(t\) noise with \(\nu=3\).}
    We use
    \[
        F(x,\xi)
        =
        \phi(x)+10^{-2}T_3,
        \qquad
        T_3\sim t_3.
    \]
    This distribution has finite variance but is not sub-Gaussian.
    It therefore falls within the finite-variance setting in
    Assumption~\ref{ass_finite_variance}.

    \item \textbf{Student-\(t\) noise with \(\nu=1\).}
    We use
    \[
        F(x,\xi)
        =
        \phi(x)+10^{-2}T_1,
        \qquad
        T_1\sim t_1.
    \]
    This is the Cauchy distribution.
    Its variance is infinite, so this experiment lies outside the
    finite-variance theory.
    We include it as a heavy-tailed stress test.

    \item    
    \textbf{Failure noise.}
\liurevision{We use
\[
    F(x,\xi)
    =
    \begin{cases}
        \phi(x),
        & \text{with probability }q,\\
        -10^{4},
        & \text{with probability }1-q,
    \end{cases}
\]
with
\[
    q=0.8.
\]

Thus, a failed evaluation returns the fixed value \(V=-10^{4}\),
as in Assumption~\ref{ass_failure}.
Each evaluation independently returns the exact value with
probability \(q=0.8>1/2\), so the majority-correctness mechanism
for the sample median applies.}
\end{itemize}

\paragraph{Sampling rules.}

For the comparisons reported below, we use the same two sampling rules
on both test sets.

\begin{itemize}

    \item \textbf{Fixed sampling.}
    We use
    \[
        N_k=30
    \]
    samples at every evaluated point.

    \item \textbf{Radius-dependent sampling.}
    We use
    \[
        N_k
        =
        \max\left\{
            1,
            \left\lfloor\frac{1}{\Delta_k}\right\rfloor
        \right\},
    \]
    where \(\Delta_k\) is the current trust-region radius.

\end{itemize}

The radius-dependent rule increases the sampling effort as the
trust-region radius decreases.
It is intended as a practical approximation to the increasing-accuracy
requirement in STORM.

The theoretical sampling requirements are considerably stronger.
For example, the concentration results in
Section~\ref{subsec:aggregation-summary}, together with the STORM
accuracy scale \revision{$\delta=\kappa_f\Delta_k^2$}
typically give \revision{$N_k=O(\Delta_k^{-4})$}
when the estimator has the usual \(O(\delta^{-2})\) accuracy
dependence.
The same \(\Delta_k^{-4}\) scaling appears in the finite-variance
sample-averaging construction in
\cite{chen2018stochastic}.

Using this theoretical rate directly would make the number of samples
grow very quickly as the trust-region radius decreases.
We therefore use the milder inverse-radius rule as a practical sampling
strategy.
The same rule has previously been tested for noisy model-based DFO in
\cite{cartis2019improving}.
We do not claim that it satisfies the theoretical STORM sampling
requirement.

\paragraph{Evaluation budgets and repeated runs.}

The two test sets use different evaluation budgets.

For a Mor\'e--Wild problem of dimension \(d\), each run is limited to \revision{$500(d+1)$}
noisy objective evaluations.
For a CUTEst problem of dimension \(d\), each run is limited to $\revision{200(d+1)}$
noisy objective evaluations.

Repeated samples at the same point are counted separately in both
budgets.
Each problem and algorithmic configuration is run independently
10 times.

The CUTEst experiments therefore use a smaller normalized evaluation
budget than the Mor\'e--Wild experiments.
For this reason, differences in the final solved fractions between the
two test sets should not be interpreted as being caused by dimension
alone.

\paragraph{Success criterion.}

The solver observes only noisy objective values.
We evaluate its progress using the noise-free objective \(\phi\).
This follows the distinction in
Section~\ref{sec:storm-background}, where \(\phi\) is the ideal
objective that we ultimately wish to minimize.

This distinction is particularly important under failure noise.
A corrupted observation can be very small even when the corresponding
point is poor for the noise-free objective.
Such an observation should not be counted as successful optimization
progress.

Let \(\mathcal P\) denote the test set under consideration.
For each \(p\in\mathcal P\), let \(\phi_p\) denote the noise-free
objective and let \(x_{p,0}\) be the initial point.
Define \revision{$\phi_p^0
    :=
    \phi_p(x_{p,0})$,}
and let \(\phi_p^{\min}\) be the best noise-free objective value found
across the methods and runs used in that test-set comparison.
A run is successful if it finds a point \(x\) satisfying \revision{$\phi_p(x)
    \leq
    \phi_p^{\min}
    +
    \tau
    \bigl(
        \phi_p^0-\phi_p^{\min}
    \bigr)$,}
where \revision{we use accuracy level $\tau=10^{-2}$ throughout.}

\paragraph{Performance profiles.}

We summarize the results using performance profiles
\cite{dolan2002benchmarking}.
Each profile is constructed separately for one test set and one sampling
rule.
The four curves in each profile therefore correspond to the four
aggregation rules.

Let \(\mathcal A\) denote the set of aggregation rules and let \(R=10\)
be the number of independent runs.
For problem \(p\), run \(r\), and aggregation rule
\(A\in\mathcal A\), let \revision{$B_{p,r}(A)$}
denote the number of noisy objective evaluations required to reach the
target.
If the target is not reached within the evaluation budget, we set \revision{$B_{p,r}(A)=\infty$.}
For each problem--run pair, define
\[
    B_{p,r}^{*}
    :=
    \min_{A\in\mathcal A} B_{p,r}(A).
\]
For \(\alpha\geq1\), the performance profile of aggregation rule \(A\)
is
\[
    \pi_A(\alpha)
    =
    \frac{1}{|\mathcal P|R}
    \left|
    \left\{
        (p,r):
        B_{p,r}(A)
        \leq
        \alpha B_{p,r}^{*}
    \right\}
    \right|.
\]

If no aggregation rule reaches the target for a particular problem--run
pair, that pair contributes zero to every profile.
Thus, the right-hand plateau also reflects the fraction of
problem--run pairs solved within the available budget.

\subsection{\revision{Results}}

\subsubsection{Gaussian noise}

\begin{figure}[htbp]
    \centering

    \begin{subfigure}[t]{0.48\textwidth}
        \centering
        \includegraphics[width=\linewidth]
        {\detokenize{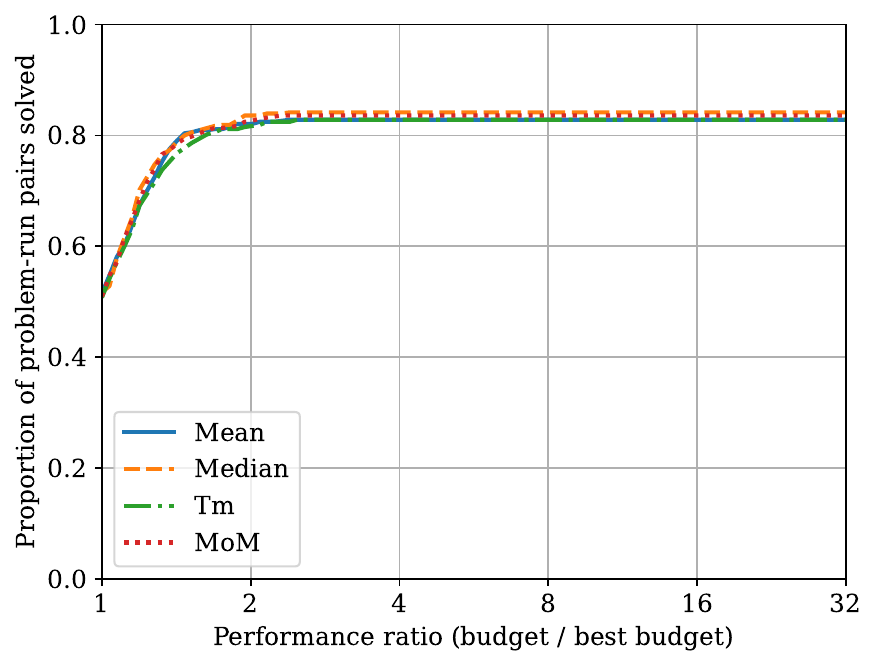}}
        \caption{Mor\'e--Wild, fixed \(N=30\).}
    \end{subfigure}
    \hfill
    \begin{subfigure}[t]{0.48\textwidth}
        \centering
        \includegraphics[width=\linewidth]
        {\detokenize{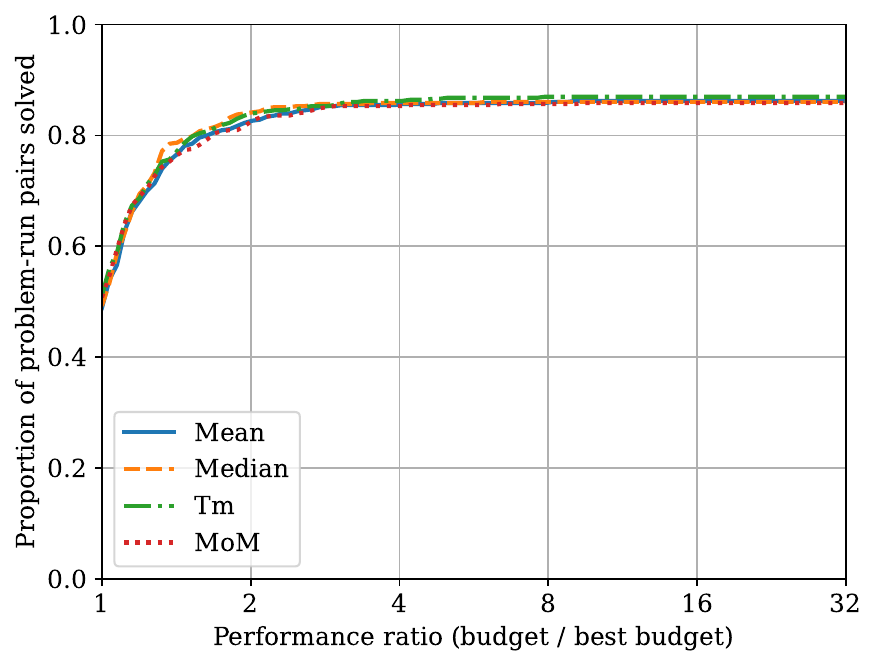}}
        \caption{Mor\'e--Wild, radius-dependent sampling.}
    \end{subfigure}

    \vspace{0.5em}

    \begin{subfigure}[t]{0.48\textwidth}
        \centering
        \includegraphics[width=\linewidth]
        {\detokenize{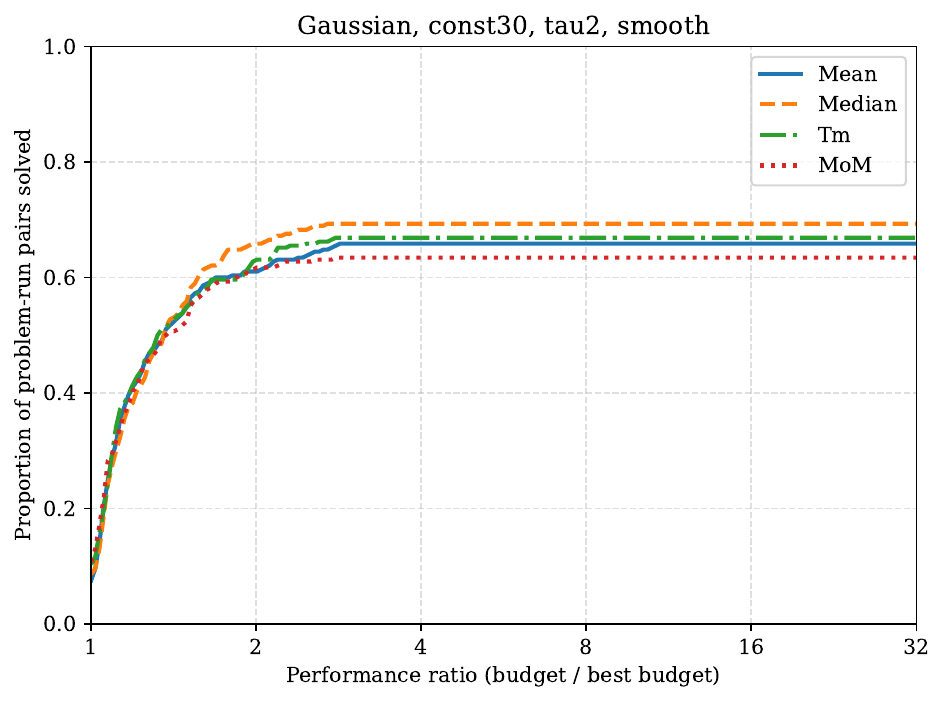}}
        \caption{CUTEst, fixed \(N=30\).}
    \end{subfigure}
    \hfill
    \begin{subfigure}[t]{0.48\textwidth}
        \centering
        \includegraphics[width=\linewidth]
        {\detokenize{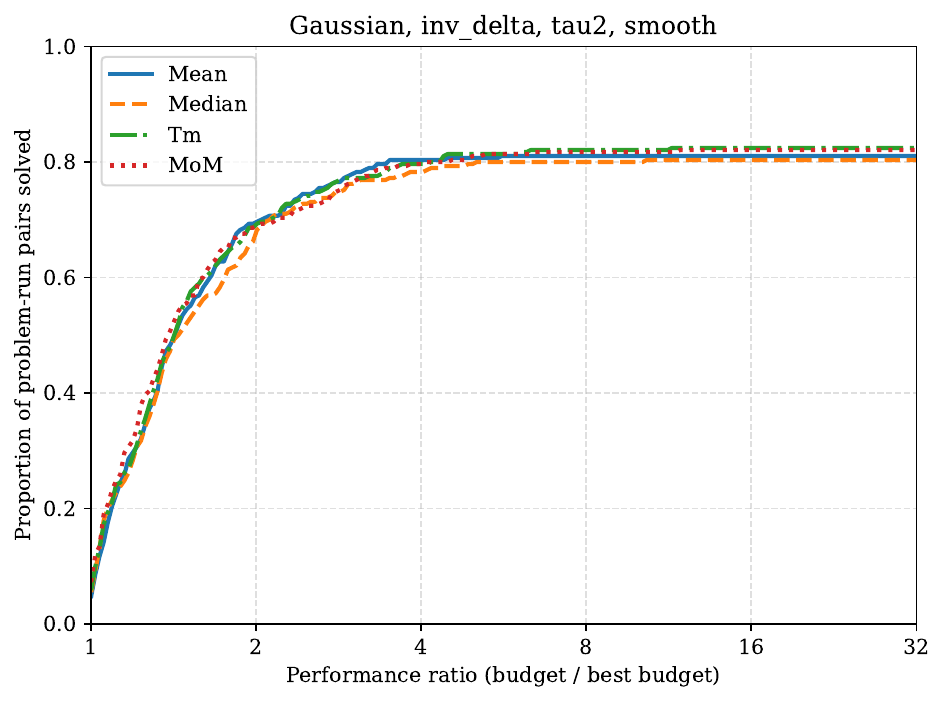}}
        \caption{CUTEst, radius-dependent sampling.}
    \end{subfigure}

    \caption{Performance profiles under Gaussian noise with standard
    deviation \(10^{-2}\) and target accuracy \(\tau=10^{-2}\).
    The top row shows the Mor\'e--Wild test set and the bottom row shows
    the CUTEst test set.}
    \label{fig:gaussian_profiles}
\end{figure}

Figure~\ref{fig:gaussian_profiles} shows a consistent pattern across
the two test sets.
The four aggregation rules perform similarly under Gaussian noise.

On the Mor\'e--Wild problems, the profiles are almost indistinguishable
for both sampling rules.
The final solved fraction is slightly above \(0.8\).

The CUTEst problems are more difficult under fixed sampling.
\liurevision{Around \(65\%\)--\(70\%\) of the problem--run pairs are solved
within the available budget.}
However, the four aggregation rules remain close.
Radius-dependent sampling substantially increases the solved fraction,
to around \(0.8\), while the differences between the aggregation rules
remain small.

This is consistent with the theory in
Section~\ref{sec:aggregation-bounds}.
Under sub-Gaussian noise, the sample mean already has exponential
concentration.
Robust aggregation therefore does not improve the order of its
high-probability sample-size guarantee.
In these experiments, the sampling schedule has a much larger effect
than the choice of aggregation rule.

\FloatBarrier

\subsubsection{Student-\texorpdfstring{\(t\)}{t} noise with
\texorpdfstring{\(\nu=3\)}{nu=3}}

\begin{figure}[htbp]
    \centering

    \begin{subfigure}[t]{0.48\textwidth}
        \centering
        \includegraphics[width=\linewidth]
        {\detokenize{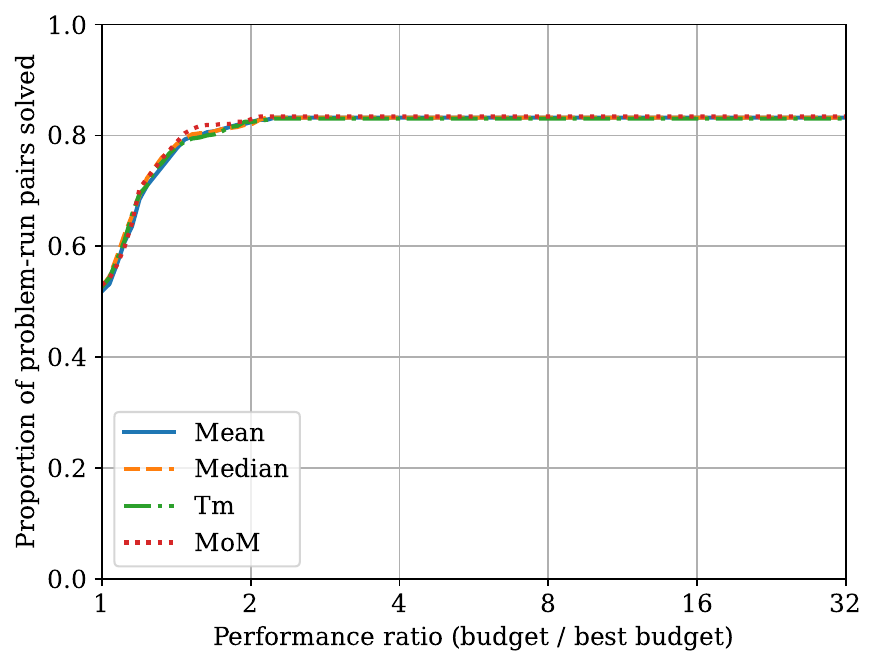}}
        \caption{Mor\'e--Wild, fixed \(N=30\).}
    \end{subfigure}
    \hfill
    \begin{subfigure}[t]{0.48\textwidth}
        \centering
        \includegraphics[width=\linewidth]
        {\detokenize{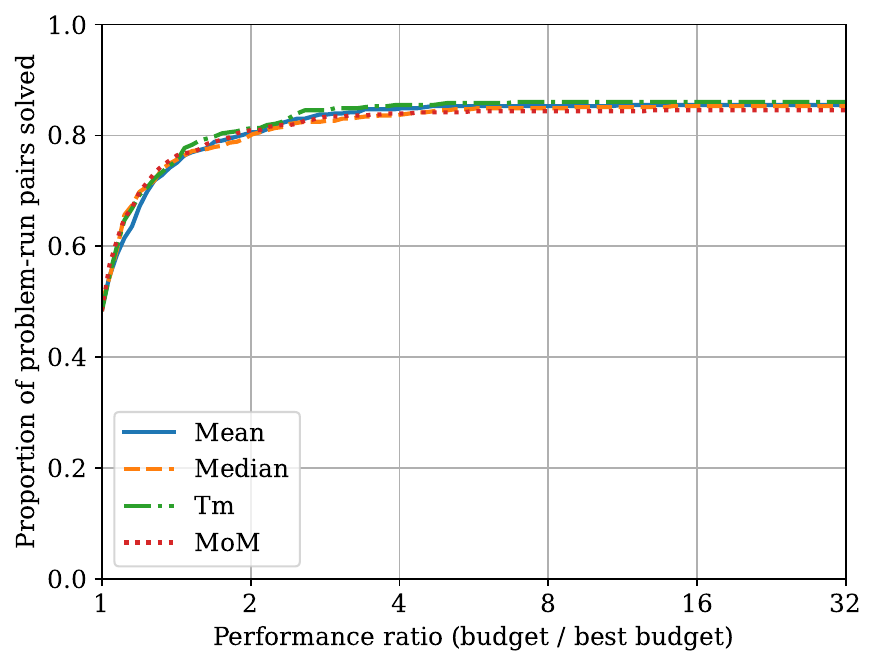}}
        \caption{Mor\'e--Wild, radius-dependent sampling.}
    \end{subfigure}

    \vspace{0.5em}

    \begin{subfigure}[t]{0.48\textwidth}
        \centering
        \includegraphics[width=\linewidth]
        {\detokenize{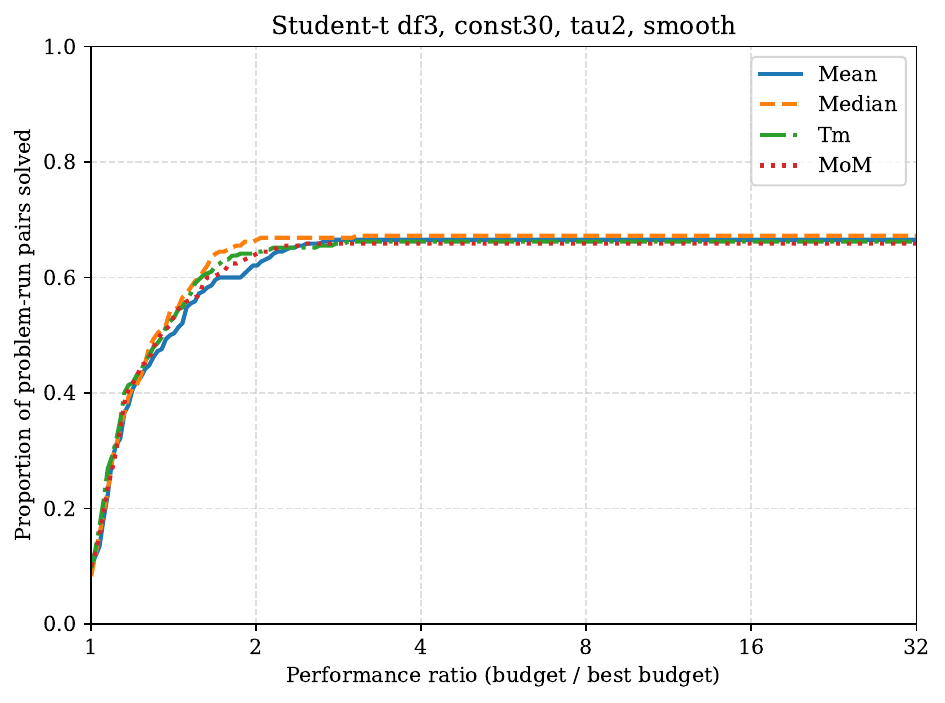}}
        \caption{CUTEst, fixed \(N=30\).}
    \end{subfigure}
    \hfill
    \begin{subfigure}[t]{0.48\textwidth}
        \centering
        \includegraphics[width=\linewidth]
        {\detokenize{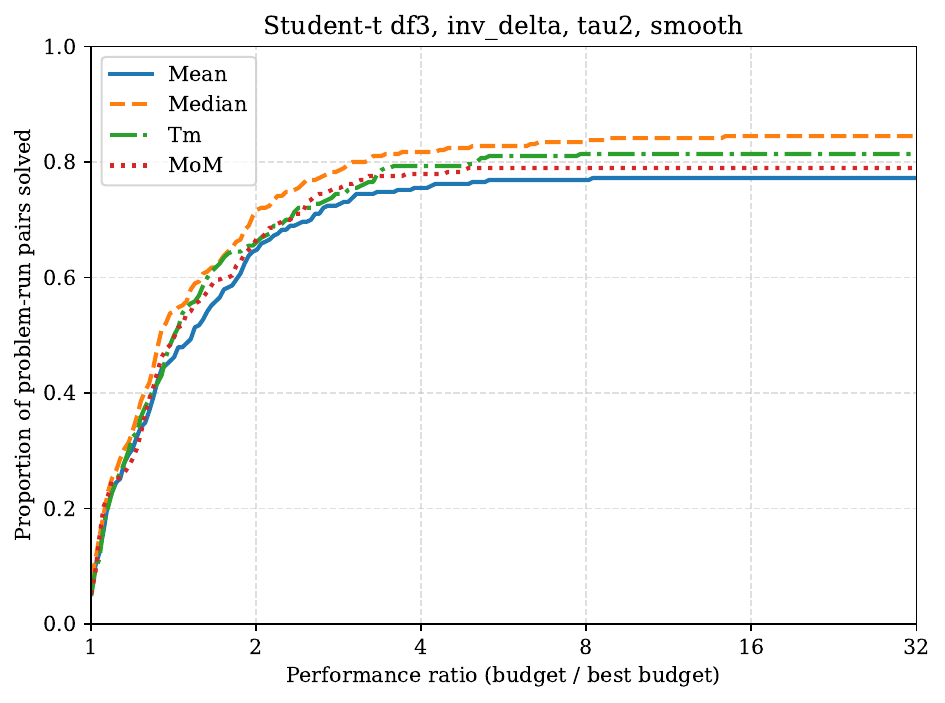}}
        \caption{CUTEst, radius-dependent sampling.}
    \end{subfigure}

    \caption{Performance profiles under Student-\(t\) noise with
    \(\nu=3\), scale \(10^{-2}\), and target accuracy
    \(\tau=10^{-2}\).
    The top row shows the Mor\'e--Wild problems and the bottom row shows
    the CUTEst problems.}
    \label{fig:student_df3_profiles}
\end{figure}

Figure~\ref{fig:student_df3_profiles} gives a similar picture under
finite-variance Student-\(t\) noise.

On the Mor\'e--Wild problems, the four profiles are again almost
indistinguishable.
This remains true for both sampling rules.

Under fixed sampling on the CUTEst problems, the four methods also
remain close.
\liurevision{Their final solved fractions are around \(65\%\)--\(70\%\).}
With radius-dependent sampling, all four methods improve substantially.
A modest separation becomes visible.
The sample median has the highest final solved fraction, while the
sample mean has the lowest.

The difference is still small compared with the heavier-tailed settings
below.
Thus, the improved high-confidence bounds for the robust estimators
under finite variance do not translate into a large solver-level
advantage in this experiment.
The CUTEst results do, however, show a modest advantage for the median
with radius-dependent sampling.

\FloatBarrier

\subsubsection{Student-\texorpdfstring{\(t\)}{t} noise with
\texorpdfstring{\(\nu=1\)}{nu=1}}

\begin{figure}[htbp]
    \centering

    \begin{subfigure}[t]{0.48\textwidth}
        \centering
        \includegraphics[width=\linewidth]
        {\detokenize{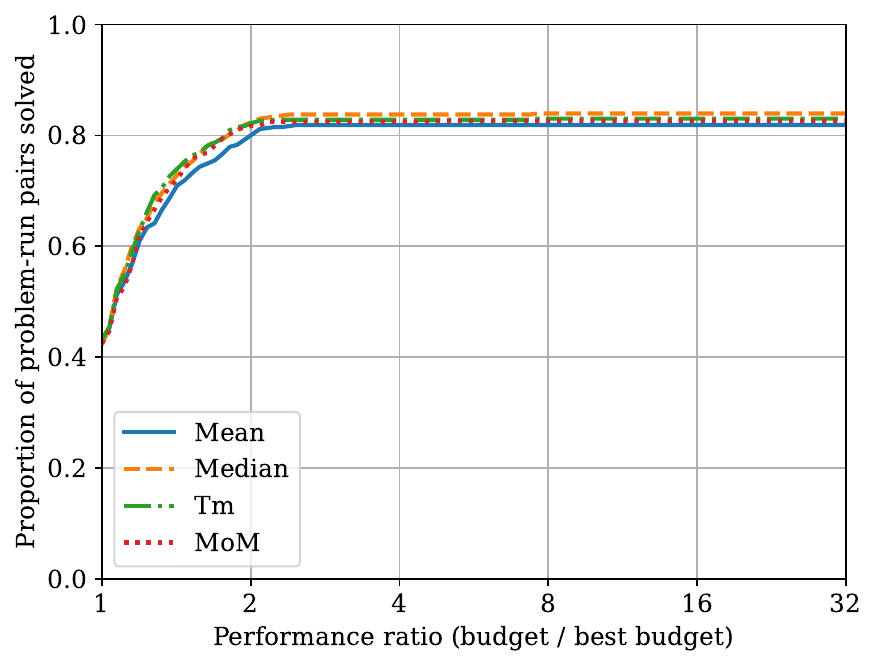}}
        \caption{Mor\'e--Wild, fixed \(N=30\).}
    \end{subfigure}
    \hfill
    \begin{subfigure}[t]{0.48\textwidth}
        \centering
        \includegraphics[width=\linewidth]
        {\detokenize{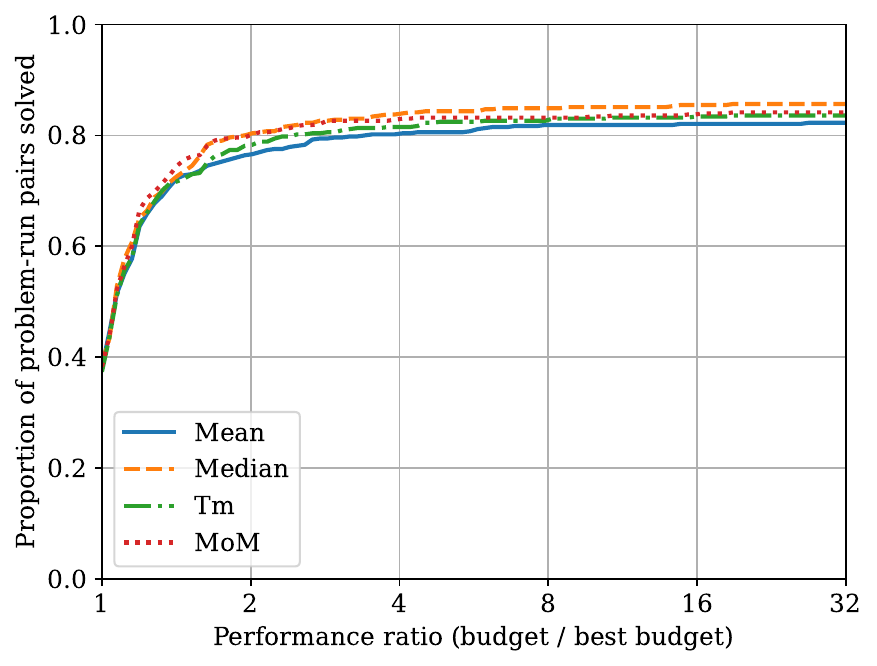}}
        \caption{Mor\'e--Wild, radius-dependent sampling.}
    \end{subfigure}

    \vspace{0.5em}

    \begin{subfigure}[t]{0.48\textwidth}
        \centering
        \includegraphics[width=\linewidth]
        {\detokenize{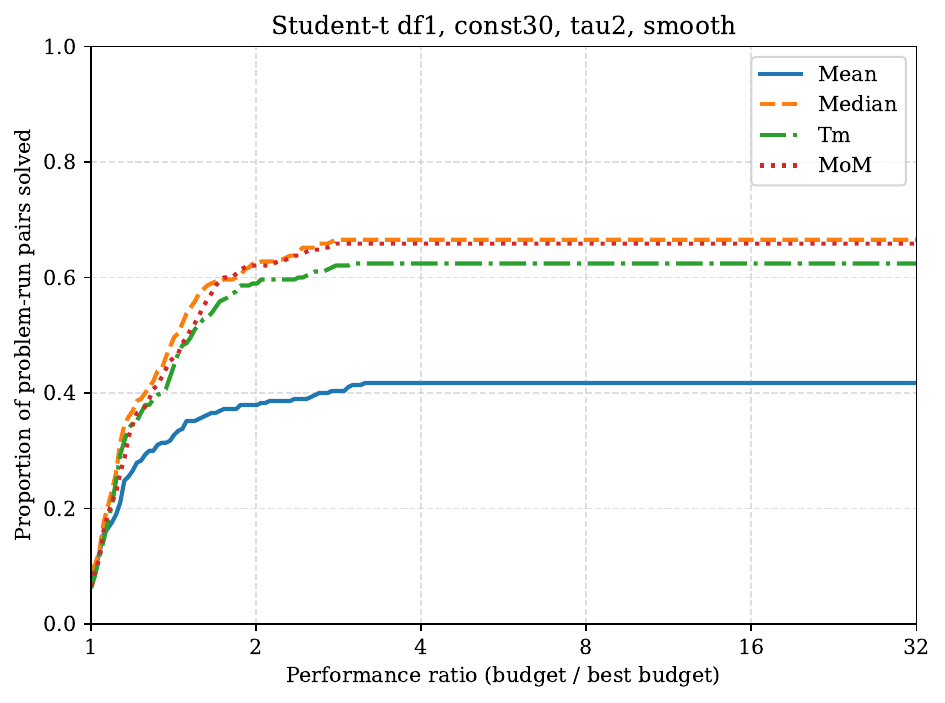}}
        \caption{CUTEst, fixed \(N=30\).}
    \end{subfigure}
    \hfill
    \begin{subfigure}[t]{0.48\textwidth}
        \centering
        \includegraphics[width=\linewidth]
        {\detokenize{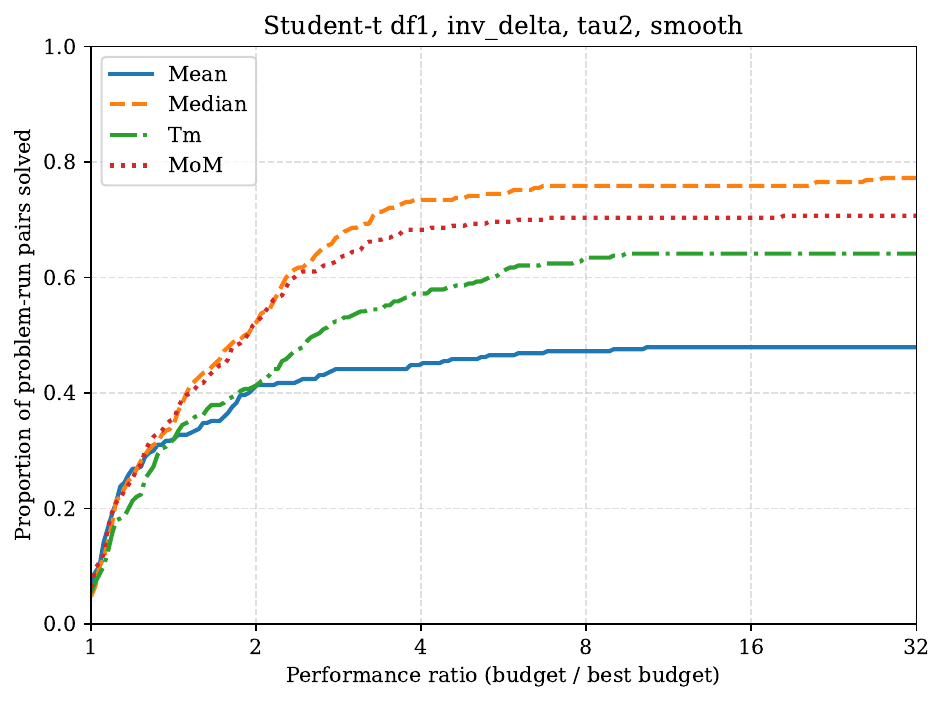}}
        \caption{CUTEst, radius-dependent sampling.}
    \end{subfigure}

    \caption{Performance profiles under Student-\(t\) noise with
    \(\nu=1\), scale \(10^{-2}\), and target accuracy
    \(\tau=10^{-2}\).
    The top row shows the Mor\'e--Wild problems and the bottom row shows
    the CUTEst problems.}
    \label{fig:student_df1_profiles}
\end{figure}

Figure~\ref{fig:student_df1_profiles} shows a much clearer difference
between the aggregation rules.

On the Mor\'e--Wild problems, the robust estimators perform only
slightly better than the sample mean at the target
\(\tau=10^{-2}\).
All four methods still solve a large fraction of the problem--run
pairs.

The CUTEst results show a much stronger separation.
\liurevision{Under fixed sampling, the sample mean solves around \(40\%\) of
the problem--run pairs.}
The robust estimators perform substantially better.
The sample median has the highest final solved fraction, followed by
median-of-means and the trimmed mean.

The separation becomes even clearer with radius-dependent sampling.
\liurevision{The sample median solves around \(80\%\) of the problem--run pairs,
and median-of-means also performs substantially better than the sample
mean.}
The sample mean reaches a final solved fraction just below \(0.5\).

The finite-variance guarantee in
Proposition~\ref{prop:mean-chebyshev} does not apply in this setting.
Cauchy observations can take extremely large values and can have a
substantial effect on an ordinary sample average.
The robust estimators reduce the influence of these observations.

The CUTEst experiment therefore strengthens the conclusion from the
Mor\'e--Wild tests.
The advantage of robust aggregation is already visible at the moderate
target \(\tau=10^{-2}\).

We do not interpret the larger separation on CUTEst as a pure dimension
effect.
The two test sets contain different objective functions and use
different evaluation budgets.
Instead, the result provides additional evidence that robust aggregation
can be important on more challenging problems under infinite-variance
noise.

\FloatBarrier

\subsubsection{Failure noise}

\begin{figure}[htbp]
    \centering

    \begin{subfigure}[t]{0.48\textwidth}
        \centering
        \includegraphics[width=\linewidth]
        {\detokenize{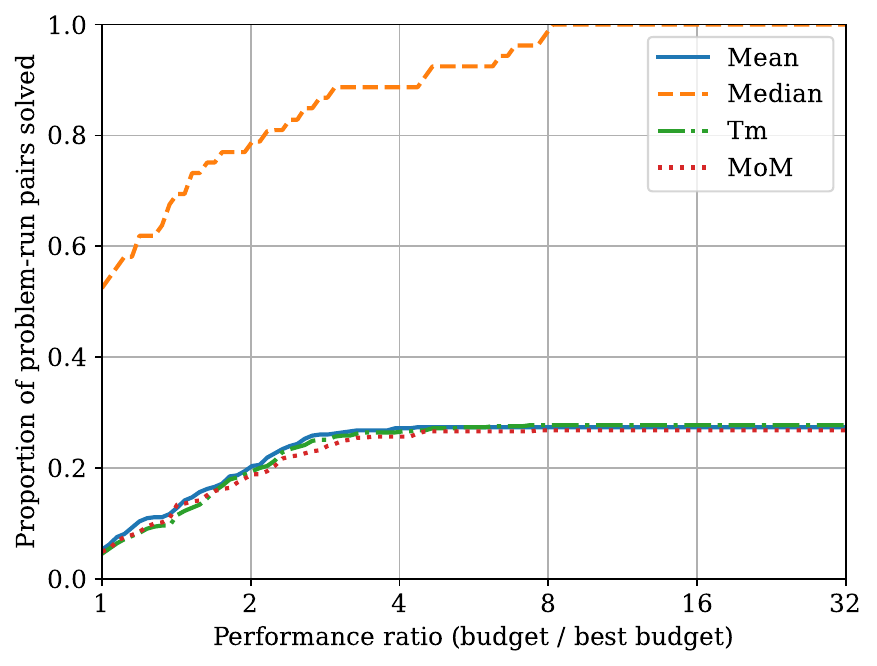}}
        \caption{Mor\'e--Wild, fixed \(N=30\).}
    \end{subfigure}
    \hfill
    \begin{subfigure}[t]{0.48\textwidth}
        \centering
        \includegraphics[width=\linewidth]
        {\detokenize{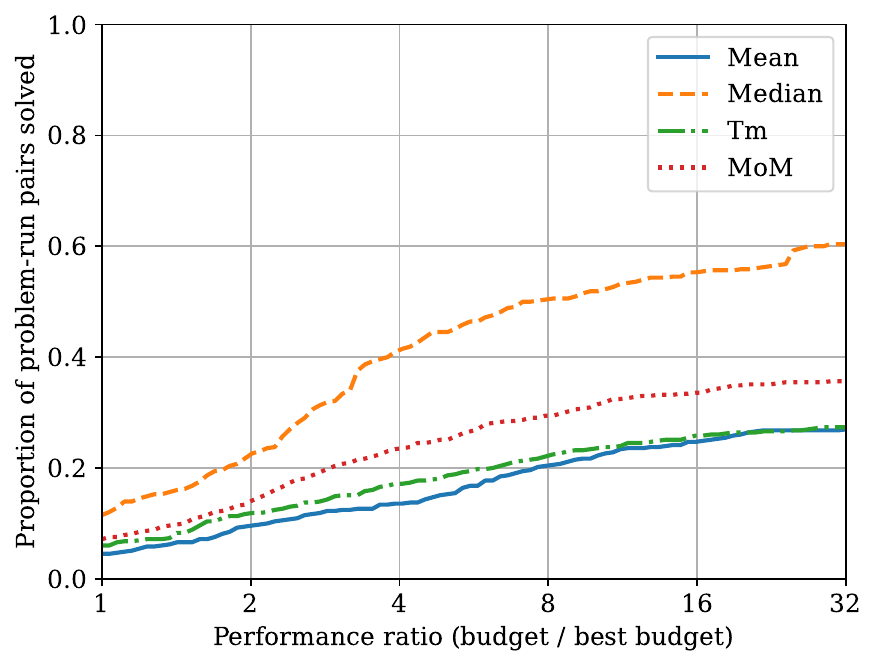}}
        \caption{Mor\'e--Wild, radius-dependent sampling.}
    \end{subfigure}

    \vspace{0.5em}

    \begin{subfigure}[t]{0.48\textwidth}
        \centering
        \includegraphics[width=\linewidth]
        {\detokenize{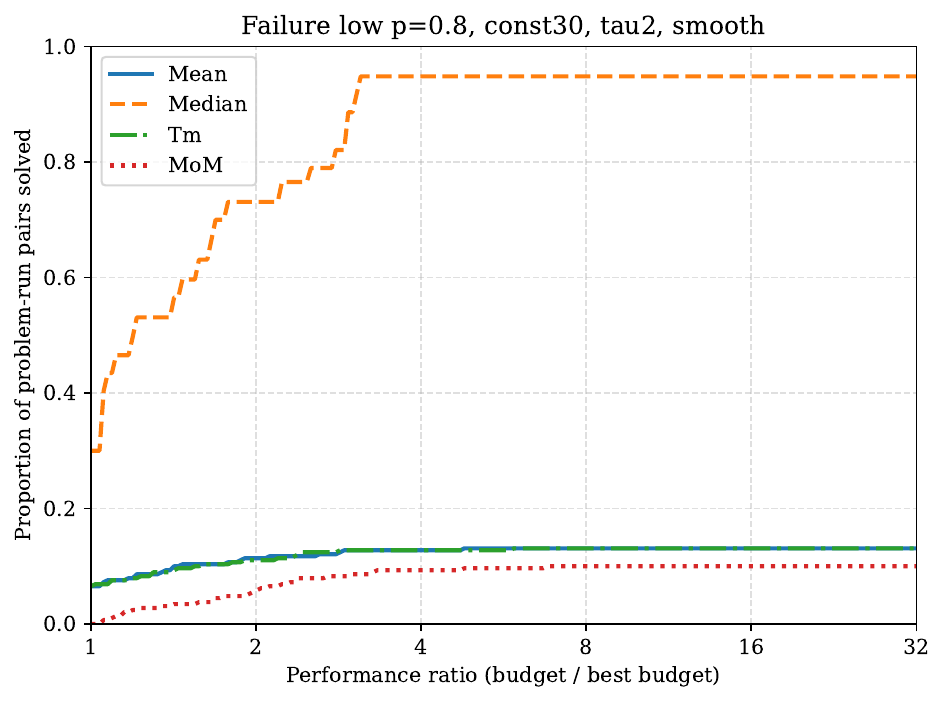}}
        \caption{CUTEst, fixed \(N=30\).}
    \end{subfigure}
    \hfill
    \begin{subfigure}[t]{0.48\textwidth}
        \centering
        \includegraphics[width=\linewidth]
        {\detokenize{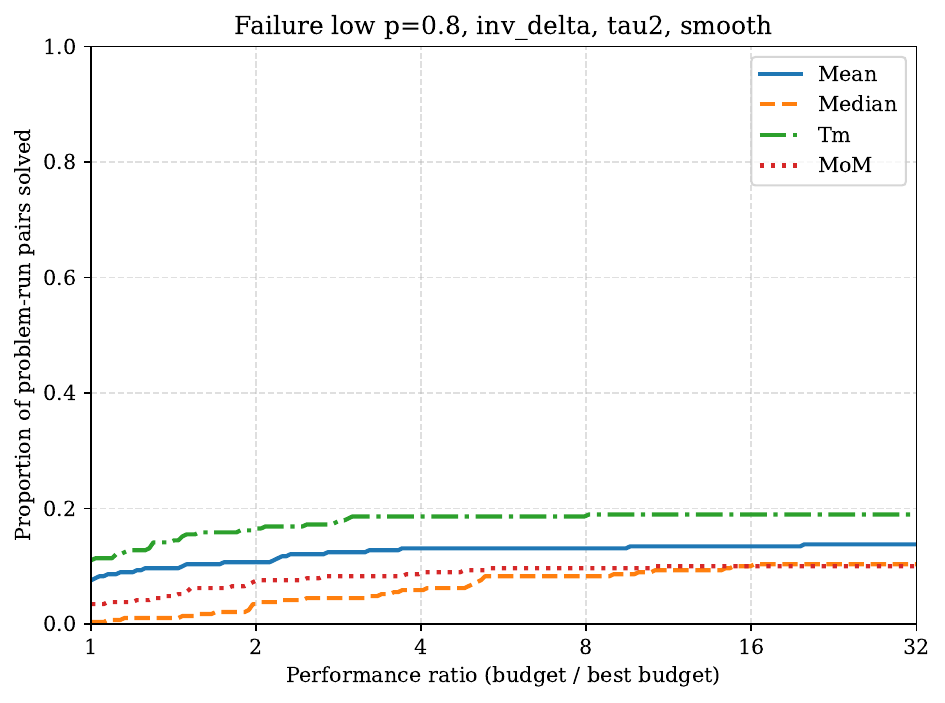}}
        \caption{CUTEst, radius-dependent sampling.}
    \end{subfigure}

\caption{\liurevision{Performance profiles under failure noise with \(q=0.8\),
fixed failed-evaluation value \(V=-10^{4}\), and target accuracy
\(\tau=10^{-2}\).
The top row shows the Mor\'e--Wild problems and the bottom row shows
the CUTEst problems.}}
\label{fig:failure_profiles}
\end{figure}

Figure~\ref{fig:failure_profiles} shows the strongest separation
between the aggregation rules.

Under fixed sampling on the Mor\'e--Wild problems, the sample median
solves almost all problem--run pairs.
The other three estimators solve only a small fraction.
The same qualitative behaviour persists on the CUTEst test set.
\liurevision{The median solves approximately \(95\%\) of the problem--run pairs,
while the other methods solve only around \(10\%\)--\(15\%\) within the available budget.}

This behaviour is closely related to
Theorem~\ref{thm:median-failure}.
Each evaluation is exact with probability
\[
    q=0.8>\frac12.
\]
With \(N=30\),
\[
    \mathbb P\left(
        \widehat f_{\mathrm{med}}(x)=\phi(x)
    \right)
    \geq
    1-\exp(-5.4)
    \approx0.995.
\]
The bound is independent of the magnitude of the corrupted
observations.
The fixed-\(N\) median can therefore exploit the majority of exact
function values effectively.

The result changes sharply under radius-dependent sampling.
On the Mor\'e--Wild problems, the median still performs better than
the other estimators, but its solved fraction drops substantially.
On the CUTEst problems, each method solves only around
\(10\%\)--\(20\%\) of the problem--run pairs, with the trimmed mean
performing best.

This behaviour shows an important interaction between aggregation and
sample allocation.
When the trust-region radius is relatively large \revision{(most commonly at the start of a solve)}, the rule
\[
    N_k
    =
    \max\left\{
        1,
        \left\lfloor\Delta_k^{-1}\right\rfloor
    \right\}
\]
may use only a small number of observations.
The sample median then has much less protection from failed
evaluations.
In particular, the high-probability majority mechanism obtained with
\(N=30\) is no longer available during these iterations.

Later in the run, when the radius becomes small, the same rule can use
many more samples.
However, these additional evaluations consume the fixed budget and may
arrive after inaccurate early models or step-acceptance decisions have
already affected the optimization path.

The exact-recovery result for the \(10\%\) trimmed mean does not apply
to the present failure configuration.
For \(q=0.8\),
Proposition~\ref{prop:tm-exact-recovery} requires
\[
    \frac{r}{N}>1-q=0.2,
\]
whereas the experiments use approximately
\[
    \frac{r}{N}=0.1.
\]
\liurevision{There is therefore no theoretical reason to expect the trimmed mean
used here to have the same exact-recovery property as the sample
median.}

\FloatBarrier

\subsection{Discussion}

The two test sets give a broadly consistent picture of when robust
aggregation is useful.

Under Gaussian noise, the four estimators behave similarly.
The sample mean already has exponential concentration in this setting,
and the robust estimators do not improve the order of its
high-confidence sample-size guarantee.
The CUTEst experiments reinforce this observation.
On that test set, changing the sampling rule has a much larger effect
than changing the aggregation rule.

The same conclusion largely holds under Student-\(t\) noise with
\(\nu=3\).
The robust estimators have better worst-case high-confidence guarantees
than the Chebyshev guarantee for the sample mean, but the solver-level
differences remain modest.
A small advantage for the sample median appears on the CUTEst problems
with radius-dependent sampling, but the four methods remain broadly
comparable.

The picture changes under Student-\(t\) noise with \(\nu=1\).
On the Mor\'e--Wild problems, robust aggregation gives only a moderate
advantage at \(\tau=10^{-2}\).
On the CUTEst test set, the difference is much larger.
The sample median, median-of-means, and trimmed mean all solve
substantially more problem--run pairs than the sample mean.
Among the four estimators, the sample median gives the strongest
performance in these experiments.

Failure noise produces the clearest separation.
With fixed \(N=30\), the sample median performs extremely well on both
test sets.
This agrees with the majority-correctness analysis in
Section~\ref{subsec:sample-median}.
The estimator can recover the noise-free value whenever successful
evaluations form a strict majority, and its guarantee does not depend
on the magnitude of the corrupted observations.

The failure-noise results also show that a robust aggregation rule
cannot be considered independently of the sampling schedule.
The median is highly effective with fixed \(N=30\), but loses most of
this advantage under radius-dependent sampling on the CUTEst test set.
The inverse-radius rule uses relatively few observations early in the
optimization, precisely when a reliable majority may be needed to
protect the local model and the step-acceptance decision from corrupted
evaluations.

Across the four noise models, the benefit of robust aggregation becomes
more visible as the observations move away from the light-tailed
setting.
The differences are small under Gaussian noise, remain modest under
finite-variance Student-\(t\) noise, become substantial under Cauchy
noise, and are most pronounced under failure noise with fixed sampling.

The comparison between Mor\'e--Wild and CUTEst should be interpreted
with some care.
The two collections contain different optimization problems, and the
CUTEst experiments use a smaller normalized evaluation budget.
The experiments therefore do not isolate problem dimension as the only
source of the observed differences.
Instead, they show that the main qualitative conclusions from the
low-dimensional experiments persist on a second and more challenging
test set.

Overall, the experiments suggest that robustness in stochastic
model-based DFO depends jointly on the aggregation rule, the noise
mechanism, and the sampling schedule.
Robust estimators offer little advantage when the sample mean is already
well concentrated.
Their benefit becomes much clearer under extreme observations and
evaluation failures, \revision{especially as the dimension of the problem increases}.
At the same time, enough samples must be allocated at the appropriate
stage of the optimization for this robustness to be useful.

\section{Conclusion}
\label{sec:conclusion}

This paper studied how repeated noisy objective evaluations should be
aggregated in stochastic model-based derivative-free optimization.
We compared the sample mean with the sample median, trimmed mean, and
median-of-means.
Under finite-variance noise, the main theoretical difference lies in the
required confidence level.
The Chebyshev guarantee for the sample mean has
\(1/(1-\beta)\) dependence, whereas the trimmed mean and
median-of-means have logarithmic dependence on \(1/(1-\beta)\), up to
problem-dependent constants.
The sample median has the same logarithmic confidence dependence when
there is sufficient probability mass around its population median.
Under failure noise, the sample median can recover the noise-free
function value with exponentially high probability when correct
evaluations occur with probability greater than \(1/2\).

The numerical experiments on the Mor\'e--Wild and CUTEst test sets show
that the benefit of robust aggregation depends strongly on the noise
model.
The four estimators perform similarly under Gaussian noise and remain
broadly comparable under finite-variance Student-\(t\) noise.
The advantage of robust aggregation becomes clearer under Cauchy noise,
especially on the CUTEst problems.
Under failure noise, the sample median performs particularly well with
fixed sampling, in agreement with the majority-correctness analysis.
However, its advantage can be lost under radius-dependent sampling.
This shows that the aggregation rule and the sampling schedule should not
be considered independently.
Overall, no estimator is uniformly best.
The sample mean remains a natural choice under light-tailed noise, while
robust estimators become more attractive when heavy tails, extreme
observations, or evaluation failures are present, \revision{especially for problems of higher dimension}.

\liurevision{These findings suggest a direction for future work:
developing adaptive sampling rules that account for the
concentration properties of different robust estimators, suitable for the ASTRO-DF framework \cite{shashaani2018astro}.}
It would also be useful to develop robust interpolation or regression
models when repeated sampling at a single point is too expensive or
cannot remove persistent errors.
More generally, the experiments show that improved pointwise
concentration does not automatically lead to the same improvement at the
solver level.
Understanding how aggregation, model construction, sampling, and
trust-region decisions interact remains an important open question.

\bibliographystyle{plain}
\bibliography{reference}

@article{amici2026stratified,
  title={Stratified Adaptive Sampling for Derivative-Free Stochastic Trust-Region Optimization},
  author={Amici, Giovanni and Shashaani, Sara and Jain, Pranav},
  journal={arXiv preprint arXiv:2604.00178},
  year={2026}
}

@article{bandeira2014convergence,
  title   = {Convergence of Trust-Region Methods Based on Probabilistic Models},
  author  = {Bandeira, Afonso S. and Scheinberg, Katya and Vicente, Luis Nunes},
  journal = {SIAM Journal on Optimization},
  volume  = {24},
  number  = {3},
  pages   = {1238--1264},
  year    = {2014},
  doi     = {10.1137/130915984}
}

@article{baraldi2025proxstorm,
  title={ProxSTORM---A Stochastic Trust-Region Algorithm for Nonsmooth Optimization},
  author={Baraldi, Robert J. and Javeed, Aurya and Kouri, Drew P. and Scheinberg, Katya},
  journal={arXiv preprint arXiv:2510.03187},
  year={2025}
}

@article{billups2013derivative,
  title   = {Derivative-Free Optimization of Expensive Functions with Computational Error Using Weighted Regression},
  author  = {Billups, Stephen C. and Larson, Jeffrey and Graf, Peter},
  journal = {SIAM Journal on Optimization},
  volume  = {23},
  number  = {1},
  pages   = {27--53},
  year    = {2013},
  publisher={SIAM},
  doi     = {10.1137/100814688}
}

@article{blanchet2019stochastic,
  title   = {Convergence Rate Analysis of a Stochastic Trust-Region Method via Supermartingales},
  author  = {Blanchet, Jose and Cartis, Coralia and Menickelly, Matt and Scheinberg, Katya},
  journal = {INFORMS Journal on Optimization},
  volume  = {1},
  number  = {2},
  pages   = {92--119},
  year    = {2019},
  publisher={INFORMS},
  doi     = {10.1287/ijoo.2019.0016}
}

@incollection{booker1998helicopter,
  author    = {Booker, Andrew J. and Dennis, Jr., J. E. and Frank, Paul D.
               and Serafini, David B. and Torczon, Virginia},
  title     = {Optimization Using Surrogate Objectives on a Helicopter Test Example},
  booktitle = {Computational Methods for Optimal Design and Control},
  editor    = {Borggaard, Jeff and Burns, John and Cliff, Eugene and Schreck, Scott},
  pages     = {49--58},
  publisher = {Birkh{\"a}user},
  year      = {1998},
  doi       = {10.1007/978-1-4612-1780-0_3}
}

@article{cao2024first,
  title   = {First- and Second-Order High Probability Complexity Bounds for Trust-Region Methods with Noisy Oracles},
  author  = {Cao, Liyuan and Berahas, Albert S. and Scheinberg, Katya},
  journal = {Mathematical Programming},
  volume  = {207},
  number  = {1},
  pages   = {55--106},
  year    = {2024},
  doi     = {10.1007/s10107-023-01999-5}
}

@article{cartis2019improving,
  author  = {Cartis, Coralia and Fiala, Jan and Marteau, Benjamin and Roberts, Lindon},
  title   = {Improving the Flexibility and Robustness of Model-Based Derivative-Free Optimization Solvers},
  journal = {ACM Transactions on Mathematical Software},
  volume  = {45},
  number  = {3},
  pages   = {32:1--32:41},
  year    = {2019},
  doi     = {10.1145/3338517}
}

@article{chang2013strong,
  author  = {Chang, Kuo-Hao and Hong, L. Jeff and Wan, Hong},
  title   = {Stochastic Trust-Region Response-Surface Method ({STRONG})---A New Response-Surface Framework for Simulation Optimization},
  journal = {INFORMS Journal on Computing},
  volume  = {25},
  number  = {2},
  pages   = {230--243},
  year    = {2013},
  doi     = {10.1287/ijoc.1120.0498}
}

@incollection{chaudhry2025complexity,
  author    = {Chaudhry, Abraar and Scheinberg, Katya},
  title     = {On Complexity of Model-Based Derivative-Free Methods},
  booktitle = {Proceedings of the International Congress of Mathematicians 2026, Volume 7: Invited Lectures (Sections 15--20)},
  pages     = {208--228},
  year      = {2026},
  publisher = {Society for Industrial and Applied Mathematics},
  doi       = {10.1137/25M1806016}
}

@article{chen2018stochastic,
  title   = {Stochastic Optimization Using a Trust-Region Method and Random Models},
  author  = {Chen, Ruobing and Menickelly, Matt and Scheinberg, Katya},
  journal = {Mathematical Programming},
  volume  = {169},
  number  = {2},
  pages   = {447--487},
  year    = {2018},
  publisher={Springer},
  doi     = {10.1007/s10107-017-1141-8}
}

@article{chen2016hidden,
  author  = {Chen, Xiaojun and Kelley, C. T.},
  title   = {Optimization with Hidden Constraints and Embedded Monte Carlo Computations},
  journal = {Optimization and Engineering},
  volume  = {17},
  number  = {1},
  pages   = {157--175},
  year    = {2016},
  doi     = {10.1007/s11081-015-9302-1}
}

@book{conn2009introduction,
  title     = {Introduction to Derivative-Free Optimization},
  author    = {Conn, Andrew R. and Scheinberg, Katya and Vicente, Luis Nunes},
  series    = {MPS-SIAM Series on Optimization},
  volume    = {8},
  publisher = {Society for Industrial and Applied Mathematics},
  address   = {Philadelphia, PA},
  year      = {2009},
  doi       = {10.1137/1.9780898718768}
}

@article{deng2009variable,
  author  = {Deng, Geng and Ferris, Michael C.},
  title   = {Variable-Number Sample-Path Optimization},
  journal = {Mathematical Programming},
  volume  = {117},
  number  = {1--2},
  pages   = {81--109},
  year    = {2009},
  doi     = {10.1007/s10107-007-0164-y}
}

@article{dolan2002benchmarking,
  title={Benchmarking optimization software with performance profiles},
  author={Dolan, Elizabeth D. and Mor{\'e}, Jorge J.},
  journal={Mathematical Programming},
  volume={91},
  pages={201--213},
  year={2002},
  doi={10.1007/s101070100263}
}

@article{fang2026heavy,
  title={High Probability Complexity Bounds of Trust-Region Stochastic Sequential Quadratic Programming with Heavy-Tailed Noise},
  author={Fang, Yuchen and Lavaei, Javad and Na, Sen},
  journal={Mathematical Programming},
  year={2026},
  publisher={Springer},
  doi={10.1007/s10107-026-02357-x}
}

@article{feng2008green,
  author  = {Feng, Wu-chun and Feng, Xizhou and Ge, Rong},
  title   = {Green Supercomputing Comes of Age},
  journal = {IT Professional},
  volume  = {10},
  number  = {1},
  pages   = {17--23},
  year    = {2008},
  doi     = {10.1109/MITP.2008.8}
}

@article{fowkes2022pycutest,
  author  = {Fowkes, Jaroslav and Roberts, Lindon and B{\H{u}}rmen, {\'A}rp{\'a}d},
  title   = {{PyCUTEst}: An Open Source Python Package of Optimization Test Problems},
  journal = {Journal of Open Source Software},
  volume  = {7},
  number  = {78},
  pages   = {4377},
  year    = {2022},
  doi     = {10.21105/joss.04377}
}

@article{gould2015cutest,
  author  = {Gould, Nicholas I. M. and Orban, Dominique and Toint, Philippe L.},
  title   = {{CUTEst}: A Constrained and Unconstrained Testing Environment with Safe Threads for Mathematical Optimization},
  journal = {Computational Optimization and Applications},
  volume  = {60},
  number  = {3},
  pages   = {545--557},
  year    = {2015},
  doi     = {10.1007/s10589-014-9687-3}
}

@article{gratton2018complexity,
  title   = {Complexity and Global Rates of Trust-Region Methods Based on Probabilistic Models},
  author  = {Gratton, Serge and Royer, Cl{\'e}ment W. and Vicente, Luis Nunes and Zhang, Zaikun},
  journal = {IMA Journal of Numerical Analysis},
  volume  = {38},
  number  = {3},
  pages   = {1579--1597},
  year    = {2018},
  doi     = {10.1093/imanum/drx043}
}

@article{ha2025adaptive,
  title={Iteration Complexity and Finite-Time Efficiency of Adaptive Sampling Trust-Region Methods for Stochastic Derivative-Free Optimization},
  author={Ha, Yunsoo and Shashaani, Sara},
  journal={IISE Transactions},
  volume={57},
  number={5},
  pages={541--555},
  year={2025},
  publisher={Taylor \& Francis},
  doi={10.1080/24725854.2024.2335513}
}

@article{ha2025complexity,
  title={Complexity of Zeroth- and First-Order Stochastic Trust-Region Algorithms},
  author={Ha, Yunsoo and Shashaani, Sara and Pasupathy, Raghu},
  journal={SIAM Journal on Optimization},
  volume={35},
  number={3},
  pages={2098--2127},
  year={2025},
  publisher={SIAM},
  doi={10.1137/24M1664484}
}

@article{ha2026adaptive,
  title={Adaptive Regularization within Trust Region Methods for Stochastic Nonconvex Optimization},
  author={Ha, Yunsoo and Shashaani, Sara and Tran-Dinh, Quoc},
  journal={arXiv preprint arXiv:2604.15457},
  year={2026}
}

@article{harding2015fault,
  author  = {Harding, Brendan and Hegland, Markus and Larson, Jay and Southern, James},
  title   = {Fault Tolerant Computation with the Sparse Grid Combination Technique},
  journal = {SIAM Journal on Scientific Computing},
  volume  = {37},
  number  = {3},
  pages   = {C331--C353},
  year    = {2015},
  doi     = {10.1137/140964448}
}

@article{howland2022parameter,
    author  = {Howland, Michael F. and Dunbar, Oliver R. A.
               and Schneider, Tapio},
    title   = {Parameter Uncertainty Quantification in an
               Idealized {GCM} With a Seasonal Cycle},
    journal = {Journal of Advances in Modeling Earth Systems},
    volume  = {14},
    number  = {3},
    pages   = {e2021MS002735},
    year    = {2022},
    doi     = {10.1029/2021MS002735}
}

@article{jin2024samplecomplexity,
  author  = {Jin, Billy and Scheinberg, Katya and Xie, Miaolan},
  title   = {Sample Complexity Analysis for Adaptive Optimization Algorithms with Stochastic Oracles},
  journal = {Mathematical Programming},
  volume  = {209},
  number  = {1--2},
  pages   = {651--679},
  year    = {2025},
  doi     = {10.1007/s10107-024-02078-z}
}

@book{kosorok2008introduction,
  author    = {Kosorok, Michael R.},
  title     = {Introduction to Empirical Processes and Semiparametric Inference},
  series    = {Springer Series in Statistics},
  publisher = {Springer},
  address   = {New York},
  year      = {2008},
  doi       = {10.1007/978-0-387-74978-5}
}

@article{larson2016stochastic,
  author  = {Larson, Jeffrey and Billups, Stephen C.},
  title   = {Stochastic Derivative-Free Optimization Using a Trust Region Framework},
  journal = {Computational Optimization and Applications},
  volume  = {64},
  number  = {3},
  pages   = {619--645},
  year    = {2016},
  doi     = {10.1007/s10589-016-9827-z}
}

@article{larson2019derivative,
  title   = {Derivative-Free Optimization Methods},
  author  = {Larson, Jeffrey and Menickelly, Matt and Wild, Stefan M.},
  journal = {Acta Numerica},
  volume  = {28},
  pages   = {287--404},
  year    = {2019},
  publisher={Cambridge University Press},
  doi     = {10.1017/S0962492919000060}
}

@article{ledigabel2024taxonomy,
  author  = {Le Digabel, S{\'e}bastien and Wild, Stefan M.},
  title   = {A Taxonomy of Constraints in Black-Box Simulation-Based Optimization},
  journal = {Optimization and Engineering},
  volume  = {25},
  number  = {2},
  pages   = {1125--1143},
  year    = {2024},
  doi     = {10.1007/s11081-023-09839-3}
}

@article{lugosi2019mean,
  title={Mean estimation and regression under heavy-tailed distributions: A survey},
  author={Lugosi, G{\'a}bor and Mendelson, Shahar},
  journal={Foundations of Computational Mathematics},
  volume={19},
  number={5},
  pages={1145--1190},
  year={2019},
  publisher={Springer},
  doi={10.1007/s10208-019-09427-x}
}

@article{ma2026projection,
  title={Projection-robust zeroth-order optimization on manifolds under heavy-tailed oracle noise},
  author={Ma, Yinpu and Li, Cunlin and Zhang, Shiyu},
  journal={Journal of King Saud University Computer and Information Sciences},
  volume={38},
  number={7},
  pages={805},
  year={2026},
  publisher={Springer},
  doi={10.1007/s44443-026-01272-7}
}

@article{massart1990tight,
  title={The Tight Constant in the {Dvoretzky--Kiefer--Wolfowitz} Inequality},
  author={Massart, Pascal},
  journal={The Annals of Probability},
  volume={18},
  number={3},
  pages={1269--1283},
  year={1990},
  doi={10.1214/aop/1176990746}
}

@article{menhorn2022trust,
  title         = {A Trust-Region Method for Derivative-Free Nonlinear Constrained Stochastic Optimization},
  author        = {Menhorn, Friedrich and Augustin, Florian and Bungartz, Hans-Joachim and Marzouk, Youssef M.},
  journal       = {arXiv preprint arXiv:1703.04156},
  year          = {2017},
  eprint        = {1703.04156},
  archivePrefix = {arXiv},
  primaryClass  = {math.OC}
}

@article{more2009benchmarking,
  title={Benchmarking derivative-free optimization algorithms},
  author={Mor{\'e}, Jorge J. and Wild, Stefan M.},
  journal={SIAM Journal on Optimization},
  volume={20},
  number={1},
  pages={172--191},
  year={2009},
  doi={10.1137/080724083}
}

@article{morvai2024countable,
  author  = {Morvai, Guszt{\'a}v and Weiss, Benjamin},
  title   = {Countable Alphabet Stationary Processes with at Least One
             Memory Word and Intermittent Estimation with Universal Rates},
  journal = {ALEA, Latin American Journal of Probability and Mathematical Statistics},
  volume  = {21},
  pages   = {1073--1099},
  year    = {2024},
  doi     = {10.30757/ALEA.v21-42}
}

@misc{oliveira2025finite,
  title={Finite-Sample Properties of the Trimmed Mean},
  author={Oliveira, Roberto I. and Orenstein, Paulo and Rico, Zoraida F.},
  year={2025},
  howpublished={arXiv preprint arXiv:2501.03694},
  eprint={2501.03694},
  archivePrefix={arXiv},
  url={https://arxiv.org/abs/2501.03694},
  note={Accepted for publication in Annales de l'Institut Henri Poincar\'e, Probabilit\'es et Statistiques}
}

@article{rinaldi2024stochastic,
  author  = {Rinaldi, Francesco and Vicente, Lu{\'i}s Nunes and Zeffiro, Damiano},
  title   = {Stochastic Trust-Region and Direct-Search Methods:
             A Weak Tail Bound Condition and Reduced Sample Sizing},
  journal = {SIAM Journal on Optimization},
  volume  = {34},
  number  = {2},
  pages   = {2067--2092},
  year    = {2024},
  doi     = {10.1137/22M1543446}
}

@misc{roberts2025interpolation,
  author        = {Roberts, Lindon},
  title         = {Introduction to Model-Based Derivative-Free Optimization},
  year          = {2025},
  howpublished  = {arXiv preprint arXiv:2510.04473},
  eprint        = {2510.04473},
  archivePrefix = {arXiv},
  primaryClass  = {math.OC},
  url           = {https://arxiv.org/abs/2510.04473}
}

@article{scheinberg2025unreliable,
  title={Stochastic Adaptive Optimization with Unreliable Inputs: A Unified Framework for High-Probability Complexity Analysis},
  author={Scheinberg, Katya and Xie, Miaolan},
  journal={arXiv preprint arXiv:2511.19411},
  year={2025}
}

@article{scriva2024challenges,
    author  = {Scriva, Giuseppe and Astrakhantsev, Nikita
               and Pilati, Sebastiano and Mazzola, Guglielmo},
    title   = {Challenges of Variational Quantum Optimization
               with Measurement Shot Noise},
    journal = {Physical Review A},
    volume  = {109},
    number  = {3},
    pages   = {032408},
    year    = {2024},
    doi     = {10.1103/PhysRevA.109.032408}
}

@article{shashaani2018astro,
  title   = {{ASTRO-DF}: A Class of Adaptive Sampling Trust-Region Algorithms for Derivative-Free Stochastic Optimization},
  author  = {Shashaani, Sara and Hashemi, Fatemeh S. and Pasupathy, Raghu},
  journal = {SIAM Journal on Optimization},
  volume  = {28},
  number  = {4},
  pages   = {3145--3176},
  year    = {2018},
  doi     = {10.1137/15M1042425}
}

@book{vershynin2018high,
  title={High-dimensional probability: An introduction with applications in data science},
  author={Vershynin, Roman},
  volume={47},
  year={2018},
  publisher={Cambridge university press}
}

@book{wainwright2019high,
  title={High-Dimensional Statistics: A Non-Asymptotic Viewpoint},
  author={Wainwright, Martin J.},
  year={2019},
  publisher={Cambridge University Press}
}

@article{ye2016inventory,
    author  = {Ye, Wenhe and You, Fengqi},
    title   = {A Computationally Efficient Simulation-Based
               Optimization Method with Region-Wise Surrogate
               Modeling for Stochastic Inventory Management of
               Supply Chains with General Network Structures},
    journal = {Computers \& Chemical Engineering},
    volume  = {87},
    pages   = {164--179},
    year    = {2016},
    doi     = {10.1016/j.compchemeng.2016.01.015}
}

@article{zhang2021concentration,
  author  = {Zhang, Huiming and Chen, Song Xi},
  title   = {Concentration Inequalities for Statistical Inference},
  journal = {Communications in Mathematical Research},
  volume  = {37},
  number  = {1},
  pages   = {1--85},
  year    = {2021},
  doi     = {10.4208/cmr.2020-0041}
}

\end{document}